\documentclass[psamsfonts]{amsart}
\usepackage{graphicx}
\usepackage{amsmath}
\usepackage{eucal}
\usepackage{amscd}
\usepackage{amssymb}
\usepackage{latexsym}
\usepackage{amsfonts}
\newtheorem{thm}{Theorem}[section]
\newtheorem{lemma}{Lemma}[section]
\newtheorem{prop}{Proposition}[section]
\newtheorem{cor}{Corollary}[section]
\newtheorem{remark}{Remark}[section]
\newtheorem{example}{Example}[section]
\numberwithin{equation}{section}

\title{Topics on Mathematical crystallography}

\author{Toshikazu Sunada}
\address{
School of Interdisciplinary Mathematical Sciences, Meiji University, 
Nakano 4-21-1, Nakano-ku, Tokyo, 164-8525 Japan
}
\email{sunada@meiji.ac.jp}
\date{}
\begin{document}

\pagestyle{myheadings}

\markright{Mathematical Crystallography}

\begin{abstract}
In July 2012 the General Assembly of the United Nations resolved that 2014 should be the International Year of Crystallography, 100 years since the award of the Nobel Prize for the discovery of X-ray diffraction by crystals. On this special occasion, we address several topics in mathematical crystallography. Especially motivated by the recent development in {\it systematic design of crystal structures} by both mathematicians and crystallographers, we discuss interesting relationships among seemingly irrelevant subjects; say, {\it standard crystal models}, {\it tight frames} in the Euclidean space, {\it rational points on  Grassmannian}, and {\it quadratic Diophantine equations}. Thus our view is quite a bit different from the traditional one in mathematical crystallography.

The central object in this article is what we call {\it crystallographic tight frames}, which are, in a loose sense, considered a generalization of {\it root systems}. We shall also pass a remark on the connections with {\it tropical geometry}, a relatively new area in mathematics, specifically with combinatorial analogues of {\it Abel-Jacobi map} and {\it Abel's theorem}.
\end{abstract}
\maketitle

\section{Introduction}

\begin{figure}[htbp]
\begin{center}
\includegraphics[width=.68\linewidth]
{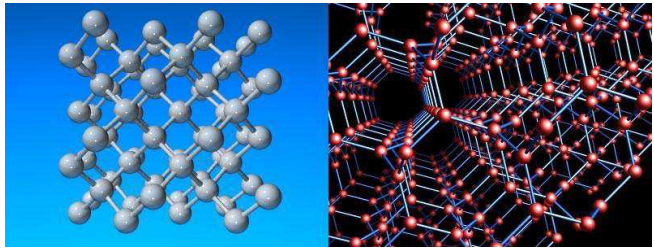}
\end{center}
\caption{Diamond and its twin}\label{fig:diamondtwin}
\end{figure}

\begin{figure}[htbp]
\begin{center}
\includegraphics[width=.92\linewidth]
{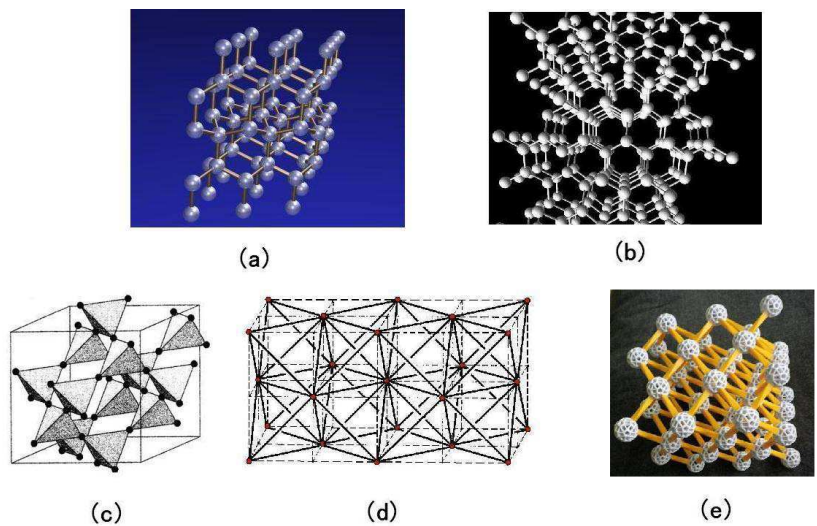}
\end{center}
\caption{(a) Lonsdaleite,~ (b) ${\rm ThSi}_{2}$ structure,~ (c) 3D kagome lattice,~ (d) Net associated with the face-centered cubic lattice, ~(e) Net associated with the body-centered cubic lattice}\label{fig:3d}
\end{figure}

It is my pastime to make various models of crystals by juggling a kit which I bought at a downtown stationer's shop. Though it is not always possible to make what I want because of the limited usage of the kit, I can still enjoy playing with it. For instance, my kit allows me to produce the model of the diamond crystal whose beauty, caused by its big symmetry, has intrigued me for some time, and motivated to look for other crystal structures, if any, with the similar symmetric property as the diamond. Actually as shown in \cite{su2} there exists the only structure that deserves to be called the {\it diamond twin}\footnote{This is what I call the {\it $K_4$~\!crystal} \cite{su2}, \cite{su5} because of the fact that it is, as a graph, the maximal abelian covering graph over the complete graph $K_4$ consisting of 4 vertices. The structure was for the first time described by Fritz Laves in 1933. Diamond and its twin are characterized by the ``strong isotropic property", the strongest one among all possible meanings of isotropy.}(Fig.~\!\ref{fig:diamondtwin}\footnote{Source of the figure of Diamond in Fig.~\!\ref{fig:diamondtwin} and Lonsdaleite in Fig.~\!\ref{fig:3d}: WebElements [http:// www.webelements.com/]. }).

In the eyes of mathematics, a crystal model as a network in space is simply a (piecewise linear) realization of an infinite-fold abelian covering graph over a finite graph. The key in this observation is that the translational action of a lattice group leaving the crystal model invariant yields a finite quotient graph\footnote{Refer to \cite{chung} for this observation made in the community of crystallographers. The {\it vector method} mentioned in this reference can be interpreted in terms of cohomology of the quotient graph (Sect.~\!\ref{sect:Standard crystal models}). Historically A. F. Wells is the crystallographer who initiated a systematic study of crystal structures as 3D networks \cite{wells1}, \cite{wells}. See \cite{nes} for some recent views on mathematical crystallography.}, and that the canonical map onto the quotient graph is a covering map whose covering transformation group is just the lattice group. This simple fact leads to the definition of {\it topological crystals} of arbitrary dimension, and can be effectively used to enumerate all topological types of crystal structures because an abelian covering graph over a finite graph $X_0$ corresponds to a subgroup of the first homology group $H_1(X_0,\mathbb{Z})$; thereby the enumeration being reduced to counting finite graphs and subgroups of their 1st homology groups \cite{su4}, \cite{su5}. Needless to say, however, there are infinitely many ways to realize the covering graph. Thus it is a natural attempt to seek a ``standard model"  having symmetries as many as possible just like Diamond and its twin.  

{\it Standard realizations} introduced in 2000 by M. Kotani and myself \cite{sk2}, \cite{sk1} in connection with asymptotic behaviors of random walks may be called standard models. Indeed the standard realization of a topological crystal $X$ has maximal symmetry in the sense that every automorphism of $X$ extends to a congruent transformation leaving the realization invariant. Moreover {\it crystal models with ``big" symmetry turn out to be the ones obtained by standard realizations} (see Theorem \ref{thm:action} in Sect.~\!\ref{sect:Standard crystal models} for the precise formulation). Figure \ref{fig:3d} illustrates several 3D examples\footnote{{\it Lonsdaleite} (named in honor of Kathleen Lonsdale) in this figure is thought of as a relative of Diamond, but is not isotropic.}. Classical 2D lattices such as the square lattice, triangular lattice, honeycomb, and kagome lattice are also standard realizations. An interesting feature following the tradition of geometry is that standard realizations are characterized by a certain minimal principle, just like the characterization of the round circle by means of the isoperimetric inequality. Furthermore this notion combined with the enumeration of topological crystals provides a useful method for a systematic design of crystal structures. Actually there is a simple algorithm for the design with which one can create a computer program to produce the CG images of two or three-dimensional crystals\footnote{Due to Hisashi Naito. Crystallographers also sought standard models; see  \cite{d-1}, \cite{d1}, \cite{eon} for instance. Some of their models are the same as ours. The algorithm {\it SYSTRE} created by Delgado--Friedrich in 2004 produces the {\it barycentric drawing}, which seems to coincide with standard realizations as far as several examples are examined. See also \cite{mc}.}.

Among all standard models, the simplest one is the {\it cubic lattice} (the 
junglegym-like figure in plain language). As a matter of fact, the cubic lattice is not much interesting as a crystal model\footnote{Sodium chloride (NaCl) crystallizes in a cubic lattice.}, but from some ``view", this lends itself to another recreation, and gives rise to an interesting mathematical issue which is linked to the standard models of general crystal structures mentioned above.

Let us look at the cubic lattice from enough remote distance. What we find out when we turn it around is that there are some specific directions toward which we may see 2D crystalline patterns (ignoring the effect of perspective). For instance, one can see the square lattice and regular triangular lattice as such crystalline patterns.

\begin{figure}[htbp]
\begin{center}
\includegraphics[width=.92\linewidth]
{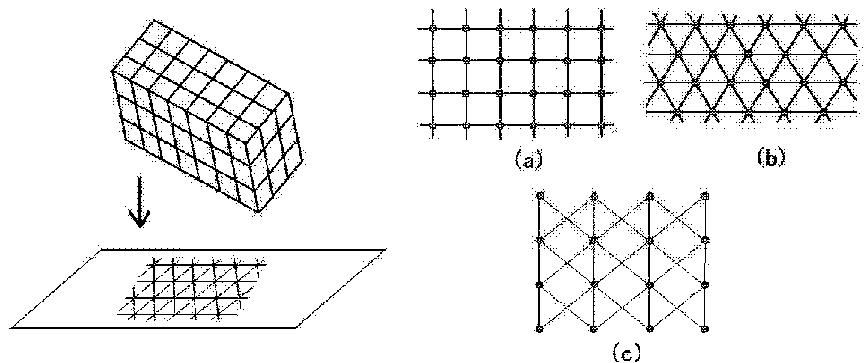}
\end{center}
\caption{Projected images of the cubic lattice}\label{fig:projection}
\end{figure}

Mathematically, we are looking at the image in the cartesian plane $\mathbb{R}^2$ of the cubic lattice placed in $\mathbb{R}^3$ by the orthogonal projection $P:(x,y,z)\mapsto (x,y)$  (see Fig.~\!\ref{fig:projection}). Here the cubic lattice is supposed to be generated by an orthonormal basis $\mathbf{f}_1,\mathbf{f}_2,\mathbf{f}_3$ of $\mathbb{R}^3$. 
Thus the set of vertices in it is
$$
\{k_1\mathbf{f}_1+k_2\mathbf{f}_2+k_3\mathbf{f}_3|~k_1,k_2,k_3\in \mathbb{Z}\}.
$$

We now put $\mathbf{v}_i=(a_i,b_i)=P(\mathbf{f}_i)$. 
Because 
\begin{equation}\label{eq001}
\langle \mathbf{x}, \mathbf{f}_1\rangle \mathbf{f}_1+\langle \mathbf{x}, \mathbf{f}_2\rangle \mathbf{f}_2+\langle \mathbf{x}, \mathbf{f}_3\rangle \mathbf{e}_3=\mathbf{x}
\quad (\mathbf{x}\in \mathbb{R}^3),
\end{equation}
and $\langle{\bf x},{\bf f}_i\rangle=\langle{\bf x},{\bf v}_i\rangle$, $P\mathbf{x}=\mathbf{x}$ for $\mathbf{x}\in \mathbb{R}^2$, 
we have, by projecting down the equality (\ref{eq001}) to the $x$-$y$ plane,
\begin{equation}\label{eq000}
\langle \mathbf{x}, \mathbf{v}_1\rangle \mathbf{v}_1+\langle \mathbf{x}, \mathbf{v}_2\rangle \mathbf{v}_2+\langle \mathbf{x}, \mathbf{v}_3\rangle \mathbf{v}_3=\mathbf{x}
\quad (\mathbf{x}\in \mathbb{R}^2).
\end{equation}

\begin{figure}[htbp]
\begin{center}
\includegraphics[width=.25\linewidth]
{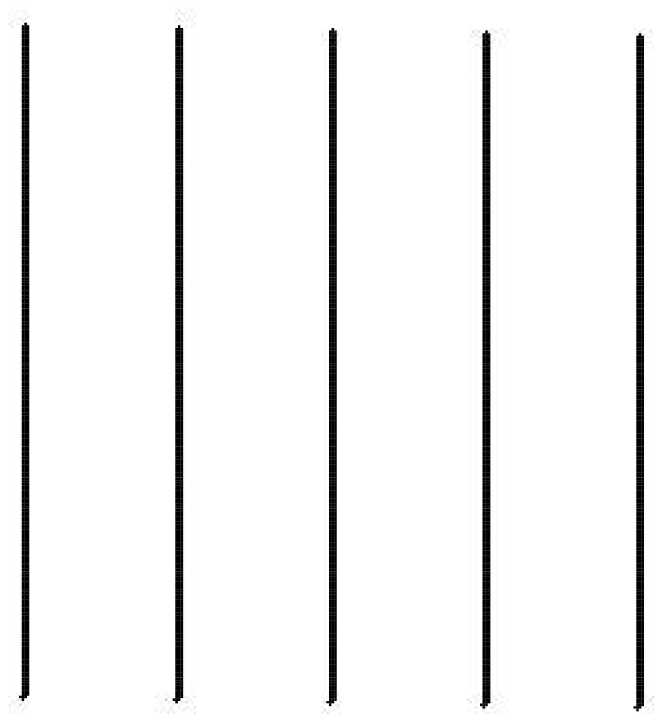}
\end{center}
\caption{}\label{fig:closure}
\end{figure}

The projected image of vertices in the cubic lattice is given by
$$
\{k_1\mathbf{v}_1+k_2\mathbf{v}_2+k_3\mathbf{v}_3|~k_1,k_2,k_3\in \mathbb{Z}\}.
$$
What we need to notice here is that the projected image does not always give a crystalline pattern. For instance, Fig.~\!\ref{fig:closure} depicts evenly spaced parallel lines expressing the closure of the projected image of vertices in the case $\mathbf{v}_1=(1,0)$, $\mathbf{v}_2=t(0,1)$, $\mathbf{v}_3=t(1,\sqrt{2})$, where we should note that $\{m+n\sqrt{2}|~m,n\in \mathbb{Z}\}$ is dense in $\mathbb{R}$ (more generally, given a positive irrational number $\alpha$, one can find infinitely many positive integers $p,q$ such that $|\alpha-q/p|<1/p^2$ (Dirichlet's theorem), from which it follows that that  $\{m+n\alpha|~m,n\in \mathbb{Z}\}$ is dense).

Actually the image of vertices in question is not discrete in almost all cases. To have a crystalline pattern, it is necessary (and sufficient) that three vectors $\mathbf{v}_1,\mathbf{v}_2,\mathbf{v}_3$ generate a lattice in $\mathbb{R}^2$, or equivalently there exist a triple of integers $(n_1,n_2,n_3)\neq (0,0,0)$ such that
\begin{equation}\label{eq01}
n_1\mathbf{v}_1+n_2\mathbf{v}_2+n_3\mathbf{v}_3=0,
\end{equation}
where one may assume without loss of generality that the greatest common divisor of $n_1,n_2,n_3$ is 1. Then the kernel of the homomorphism $\rho:\mathbb{Z}^3\longrightarrow \mathbb{R}^2$ defined by $\rho(k_1,k_2,k_3)=k_1\mathbf{v}_1+k_2\mathbf{v}_2+k_3\mathbf{v}_3$ coincides with $H:=\mathbb{Z}(n_1,n_2,n_3)$. Going back to Fig.~\!\ref{fig:projection}, we observe that the square lattice (a), regular triangular lattice (b) and the lattice (c) correspond to 
$
(n_0,n_1,n_2)=(1,0,0)$,~
$(n_0,n_1,n_2)=(1,1,1)$,~
$(n_0,n_1,n_2)=(1,1,2)$,
respectively.
\medskip

This expository article, including a few new results, is thought of as a continuation of my book \cite{su5} published in 2012. 
The purpose is, starting from the above elementary observations, to share a link to a few mathematical subjects, say {\it tight frames} in the Euclidean spaces, {\it rational points on Grassmannians}, and {\it quadratic Diophantine equations}. Those subjects are not things of novelty (for instance, tight frames appear in various guises in practical sciences), but turn out to be closely connected with each other in an interesting way.

\begin{figure}[htbp]
\begin{center}
\includegraphics[width=.5\linewidth]
{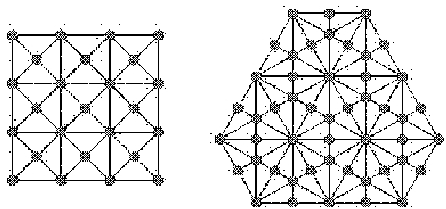}
\end{center}
\caption{Examples of 2D Coxeter complexes}\label{fig:coxeter}
\end{figure}

The protagonist is {\it crystallographic tight frames} introduced in Sect.~\!3. This notion generalizes the above-mentioned situation, and is closely related to a systematic design of crystal structures through the notion of standard realizations (Sect.~\!\ref{sect:Standard crystal models}). Further this in a special case is regarded as a generalization of {\it root systems} whose origin is in the work of W. Killings, E. Cartan and H. Weyl on Lie groups. Actually irreducible root systems yield highly symmetrical crystallographic tight frames. It should be pointed out that root systems pertain to {\it Euclidean Coxeter complexes} (cf.~\!\cite{brown}), very remarkable triangulations of the Euclidean space, named after H. S. M. Coxeter (see Fig.~\!\ref{fig:coxeter} for 2D examples which correspond to the root systems $B_2$ and $G_2$). A remarkable fact is that the 1-skeleton of a Coxeter complex is the standard realization of a crystal structure (Sect.~\!\ref{sect:Standard crystal models}).

As is well known, root systems are completely classified by means of {\it Dynkin diagrams}. On the other hand, as will be described in Sect.~\!4, 
similarity classes of crystallographic tight frames are parameterized by rational points on Grassmannians (this is by no means surprising if we rephrase the above observation as ``the projection of the cubic lattice towards a {\it rational direction} gives rise to a crystal pattern"). Certain quadratic Diophantine equations show up when we explicitly associate crystallographic tight frames with rational points. A noteworthy situation occurs in the 2-dimensional case especially; we may parameterize the (oriented) congruence classes by ``rational points" on a certain complex projective quadric.
A rational point we mean here is a point in a complex projective space each of whose homogeneous coordinate is represented by a number in an {\it imaginary quadratic field}. 
In Sect.~\!6, we explain a relationship with {\it tropical geometry}, skeletonized version of algebraic geometry, especially with discrete Abel-Jacobi maps. This unexpected link is brought about via crystallographic tight frames associated with finite graphs. The final section presents a link between discrete Abel-Jacobi maps and standard realizations.

Not surprisingly, the materials in this paper have something to do with, not only  the legacy of Eukleides--Archimedes--Kepler's achievements in polyhedral geometry\footnote{Legend has it that Pythagoras derived the notion of regular polyhedra from the shape of a crystal.}, but also with the {\it geometry of numbers} and the theory of {\it quadratic forms} because of the crucial role played by lattice groups in crystallography. 

\medskip

Before leaving the introduction, let us fix a few notations used throughout. We express 
a matrix $A$ by $(a_{ij})$ for simplicity when the $(i,j)^{\rm th}$ entry of $A$ is $a_{ij}$. The transpose of $A=(a_{ij})$ expressed by ${}^t\!A$ is the matrix whose $(i,j)^{\rm th}$ entry is $a_{ji}$. The {\it trace} of a square matrix $A$, the sum of the diagonal entries of $A$, is denoted by ${\rm tr}~\!A$. The {\it determinant} of $A$ is denoted by $\det A$.

Given a field $K$, we think of $K^d$ as a vector space over $K$ consisting of column vectors ${\bf x}={}^t(x_1,\ldots,x_d)$ with $x_i\in K$. The subspace spanned by vectors ${\bf x}_1,\ldots,{\bf x}_k\in K^d$ is expressed by $\langle {\bf x}_1,\ldots,{\bf x}_k\rangle_{K}$. We denote by $M_{m,n}(K)$ the set of all $m\times n$ matrices whose entries belong to $K$. We also use the notation $M_n(K)$ for $M_{n,n}(K)$. 
The identity matrix $(\delta_{ij})$ in $M_{n}(K)$ is denoted by $I_n$, or simply $I$. 

A matrix $A\in M_{m,n}(K)$ is identified with the linear operator of $K^n$ into $K^m$ given by ${\bf x}\mapsto A{\bf x}$. 
For a linear operator $T:K^n\longrightarrow K^m$, the {\it kernel} of $T$ is written as ${\rm ker}~\!T$. The {\it image} of $T$ is denoted by ${\rm Image}~\!T$. 
When the matrix $A\in M_{m,n}(K)$ consists of column vectors ${\bf a}_i$ ($i=1,\ldots,n$), we write $A=({\bf a}_1,\ldots,{\bf a}_n)$.

The (standard) inner product $\langle{\bf x},{\bf y} \rangle$ of two vectors ${\bf x}={}^t(x_1,\ldots,x_n)$, ${\bf y}={}^t(y_1,\ldots,y_n)$ in $\mathbb{R}^n$ is ${}^t{\bf x}{\bf y}=x_1y_1+\cdots+x_ny_n$. The norm $\|{\bf x}\|$ is $\langle {\bf x},{\bf x}\rangle^{1/2}$. For a subspace $W\subset \mathbb{R}^n$, the orthogonal complement $\big\{{\bf x}\in \mathbb{R}^n|~\langle{\bf x},{\bf y} \rangle=0 ~({\bf y}\in W)\big\}$ is denoted by $W^{\perp}$. 

\bigskip

\noindent{\bf Acknowledgement}.~~ I would like to thank Peter Kuchment for a valuable hint about the relation between standard realizations and tight frames. I also thank Hisashi Naito and my daughter Kayo for producing several figures.

\section{Tight frames}\label{sec:tight}

Property (\ref{eq000}) tells that $\mathbf{v}_1,\mathbf{v}_2,\mathbf{v}_3$ form a  {\it tight frame} of $\mathbb{R}^2$, a terminology originally used in wavelet analysis. The basic philosophy of tight frames is that representations such as (\ref{eq000}) are similar to an orthogonal expansion when considering an infinite dimensional Hilbert space such as $L^2(\mathbb{R}^d)$, and that one may have more freedom in choosing the ${\bf v}_j$ to have desirable properties such as certain smoothness and small support properties that may be impossible were they to be orthogonal (see \cite{duf} for a pioneer work). Tight frames are intimately related to rank-one quantum measurements \cite{eldar}.  In the finite dimensional case they are seen in the study of packet-based communication systems (refer to \cite{goyal} for instance), and also show up as spherical 2-designs in combinatorics \cite{waldron}. In this paper, we shall give a completely different view to tight frames. Our development is guided by the idea indicated in Introduction.

Let us start with some fundamental facts on tight frames which are more or less known (cf.~\!\cite{vale}). Only some rudiments of linear algebra is required to read this section.

In general, a sequence of $N$ vectors $\mathcal{S}=\{\mathbf{v}_i\}_{i=1}^N$ in $\mathbb{R}^d$ is said to be a $d$-dimensional frame of size $N$ or simply {\it frame} if it generates $\mathbb{R}^d$. In this definition, some of ${\bf v}_i$ allow to be zero or parallel. Given a frame, we may associate a linear operator (called the {\it frame operator}) $S=S_{\mathcal{S}}:\mathbb{R}^d\longrightarrow \mathbb{R}^d$ by setting
$$
S({\bf x})=\sum_{i=1}^N \langle{\bf x},{\bf v}_i\rangle {\bf v}_i,
$$
which is symmetric and positive. The matrix for $S$ is given by  
$$
S=({\bf v}_1,\ldots,{\bf v}_N)~\!{}^t({\bf v}_1,\ldots,{\bf v}_N)={\bf v}_1{}^t{\bf v}_1+\cdots+{\bf v}_N{}^t{\bf v}_N,
$$
and hence
\begin{equation}\label{eq:trace}
{\rm tr}~\!S
=\sum_{i=1}^N\|{\bf v}_i\|^2.
\end{equation}

A frame $\mathcal{S}$ is said to be $\alpha$-{\it tight} (or simply tight) if $S=\alpha I_d$ with a positive $\alpha$, i.e.,
$$
\sum_{i=1}^N \langle{\bf x},{\bf v}_i\rangle {\bf v}_i=\alpha{\bf x}\quad ({\bf x}\in \mathbb{R}^d).
$$ 
In view of (\ref{eq:trace}), if $\mathcal{S}$ is 1-tight, then 
$$
\sum_{i=1}^N\|{\bf v}_i\|^2=d.
$$
Tightness (resp. 1-tightness) is obviously preserved by similar transformations and permutations of subscripts $i$ in ${\bf v}_i$ (resp. by orthogonal transformations). Here two frames $\mathcal{S}_1=\{{\bf u}_i\}_{i=1}^N$ and $\mathcal{S}_2=\{{\bf v}_i\}_{i=1}^N$ are said to be {\it similar} if there exists an orthogonal transformation $U$ of $\mathbb{R}^d$ and a positive number $\lambda$ such that ${\bf u}_i=\lambda U({\bf v}_i)$ ($i=1,\ldots,N$). If $\lambda=1$ in this relation, $\mathcal{S}_1$ and $\mathcal{S}_2$ are said to be {\it congruent}.
\medskip

A tight frame appears in the following situation. 

\begin{prop}\label{thm:tightirre}
Suppose that a finite group $G$ acts on $\mathbb{R}^d$ as orthogonal transformations, and let $\mathcal{S}=\{{\bf v}_i\}_{i=1}^N$ be a frame in $\mathbb{R}^d$ which is invariant under the $G$-action (precisely speaking, for any $g\in G$, there exists a permutation $\sigma$ of $\{1,\ldots,N\}$ such that $g{\bf v}_i={\bf v}_{\sigma i}$). If the $G$-action on $\mathbb{R}^d$ is {\it irreducible}, 
then $\mathcal{S}$ is a tight frame satisfying $\displaystyle\sum_{i=1}^N{\bf v}_i={\bf 0}$. 

\end{prop}

To check this, we note that the $G$-action commutes with the frame operator $S$. Looking at eigenspaces of $S$, we conclude that $S=\alpha I_d$ for some positive scalar $\alpha$, as claimed. For the claim $\displaystyle\sum_{i=1}^N{\bf v}_i={\bf 0}$, we only have to notice that the left hand side is a $G$-invariant vector. 

\medskip

\begin{figure}[htbp]
\begin{center}
\includegraphics[width=.32\linewidth]
{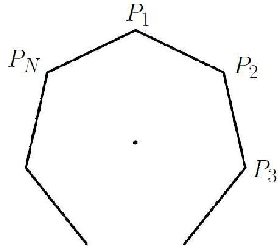}
\end{center}
\caption{A regular polygon}\label{fig:polygon}
\end{figure}

Using Proposition \ref{thm:tightirre}, one can prove that for points $P_1,\ldots, P_N$ ($N\geq 3$) in the plane $\mathbb{R}^2$ forming a $N$-regular polygon with barycenter $O$ (Fig.~\!\ref{fig:polygon}), the vectors ${\bf v}_1=\overrightarrow{OP_1},\ldots,{\bf v}_N=\overrightarrow{OP_N}$ yield a tight frame (of course, one can prove this by a direct computation). By the same reasoning, five Platonic solids (regular convex polyhedra) and thirteen Archimedean solids (semi-regular polyhedra) yield  tight frames of $\mathbb{R}^3$ (cf.~\!\cite{cromwell}).

Tight frames have a similar nature as orthonomal basis as seen in the following proposition.

\begin{prop} The following three conditions are equivalent to the 1-tightness of a frame $\mathcal{S}=\{{\bf v}_i\}_{i=1}^N$, respectively. 

{\rm (1)}~ \quad 
$\displaystyle\sum_{i=1}^N\langle{\bf v}_i,{\bf x}\rangle^2=\|{\bf x}\|^2 \quad ({\bf x}\in \mathbb{R}^d)$,

{\rm (2)}~ \quad $\displaystyle\sum_{i=1}^N\langle{\bf v}_i,{\bf x}\rangle \langle{\bf v}_i,{\bf y}\rangle=\langle {\bf x},{\bf y}\rangle\quad ({\bf x},{\bf y}\in \mathbb{R}^d)
$,

{\rm (3)}~ \quad
$\displaystyle\sum_{j=1}^N\langle T\mathbf{v}_j,\mathbf{v}_j\rangle= {\rm tr}~\!T
$
for any linear transformation $T$ of $\mathbb{R}^d$.  
\end{prop}

\noindent {\it Proof}.~ We only show that  (3) implies the 1-tightness because the other claims are easy to check. Consider the operator $T$ defined by $T({\bf x})=\langle {\bf x},{\bf y}\rangle{\bf y}$. Evidently ${\rm tr}~\!T=\|{\bf y}\|^2$ and 
$$
\sum_{j=1}^N\langle T\mathbf{v}_j,\mathbf{v}_j\rangle=\sum_{j=1}^N\langle
{\bf y},{\bf v}_j \rangle^2.
$$
Applying the equality (3), we have 
$\displaystyle\sum_{j=1}^N\langle
{\bf y},{\bf v}_j \rangle^2=\|{\bf y}\|^2$. 
Thus by (2), $\mathcal{S}$ is 1-tight.  \hfill $\Box$

\medskip

Recall that the 1-tight frame $\{\mathbf{v}_1,\mathbf{v}_2, \mathbf{v}_3\}$ mentioned in Introduction was obtained as the projected image of an orthonormal basis of $\mathbb{R}^3$. This is true for general  1-tight frames; that is, any 1-tight frame $\mathcal{S}=\{\mathbf{v}_i\}_{i=1}^N$ in $\mathbb{R}^d$ is a projected image of an orthonormal basis of $\mathbb{R}^N$. To see this, we shall introduce an auxiliary operator (matrix). 
Given a frame $\mathcal{S}=\{{\bf v}_i\}_{i=1}^N$, define the {\it frame projection} 
$
P=P_{\mathcal{S}}:\mathbb{R}^N\longrightarrow \mathbb{R}^d
$ 
by 
$$
P\big({}^t(x_1,\ldots,x_N)\big)=x_1{\bf v}_1+\cdots+x_N{\bf v}_N,
$$ 
i.e., $P$ is the linear operator characterized by $P({\bf f}_i)={\bf v}_i$ ($i=1,\ldots,N$) where $\{{\bf f}_i\}_{i=1}^N$ is the fundamental basis of $\mathbb{R}^N$. The matrix for $P$ is nothing but $({\bf v}_1,\ldots,{\bf v}_N)$. 
Therefore we have $S=P~\!{}^t\!P$. On the other hand, the matrix for ${}^t\!PP$ is the {\it Gramm matrix} $G_{\mathcal{S}}=\big(\langle {\bf v}_i,{\bf v}_j\rangle\big)$.

The next proposition is considered a special case of {\it Naimark's dilation theorem} in the theory of quantum measurements.

\begin{prop}\label{prop:nnn} The following four conditions are equivalent:

{\rm (1)}~ $\mathcal{S}$ is 1-tight.

{\rm (2)}~ $P=P_{\mathcal{S}}$ satisfies $P{}^t\!P=I$.

{\rm (3)}~ $P$ is orthogonal in the sense that 
the restriction $P\big|\big({\rm ker}~\!P\big)^{\perp}:\big({\rm ker}~\!P\big)^{\perp}$ $\longrightarrow \mathbb{R}^d$ is an isometry (i.e., it preserves the inner products).

{\rm (4)}~ 
${}^t\!PP:\mathbb{R}^N\longrightarrow\mathbb{R}^N$ is a orthogonal projection, or equivalently the Gramm matrix satisfies $G_{\mathcal{S}}^2=G_{\mathcal{S}}$.

\end{prop}

\noindent{\it Proof}.~(1) $\Leftrightarrow$ (2) is obvious. To show 
(1) $\Leftrightarrow$ (3), suppose that $\mathcal{S}$ is 1-tight. 
Since  $P~\!{}^t\!P=I$, we find
$$
\langle{}^t\!P({\bf x}),{}^t\!P({\bf y})\rangle=\langle P~\!{}^t\!P({\bf x}),{\bf y}\rangle=\langle {\bf x},{\bf y}\rangle
$$
and
$$
\langle P~\!{}^t\!P({\bf x}),P~\!{}^t\!P({\bf y})\rangle=\langle{\bf x},{\bf y} \rangle.
$$
Hence 
$$
\langle P~\!{}^t\!P({\bf x}),P~\!{}^t\!P({\bf y})\rangle=\langle{}^t\!P({\bf x}),{}^t\!P({\bf y})\rangle.
$$
Since $\big({\rm ker}~\!P\big)^{\perp}={\rm Image}~\!{}^t\!P$, we conclude that $P\big|\big({\rm ker}~\!P\big)^{\perp}:\big({\rm ker}~\!P\big)^{\perp}\longrightarrow \mathbb{R}^d$ is an isometry.

Next suppose that $P\big|\big({\rm ker}~\!P\big)^{\perp}:\big({\rm ker}~\!P\big)^{\perp}\longrightarrow \mathbb{R}^d$ is an isometry. Again using $\big({\rm ker}~\!P\big)^{\perp}={\rm Image}~\!{}^tP$, we have
$$
\langle P~\!{}^t\!P({\bf x}),P~\!{}^t\!P({\bf y})\rangle=\langle{}^t\!P({\bf x}),{}^t\!P({\bf y})\rangle,
$$
or equivalently
$$
\langle S^2({\bf x}),{\bf y}\rangle=\langle S({\bf x}),{\bf y}\rangle.
$$
This implies $S^2=S$, and hence $S=I$. 

(1) $\Leftrightarrow$ (4): If $\mathcal{S}$ is 1-tight, then 
$$
({}^t\!PP)({}^t\!PP)={}^t\!P(P{}^t\!P)P={}^t\!PP.
$$
Conversely if $({}^t\!PP)^2={}^t\!PP$, then $P{}^t\!PP=P$ since ${}^tP$ is injective. But $P$ is surjective, so $P{}^t\!P=I$. 
\hfill $\Box$
\medskip

For a frame $\mathcal{S}$ of $\mathbb{R}^d$, the $(N-d)$-dimensional subspace $W(\mathcal{S})={\rm ker}~\!P_{\mathcal{S}}$ of $\mathbb{R}^N$ is called the {\it vanishing subspace} for $\mathcal{S}$, which obviously depends only on the similarity class of $\mathcal{S}$. 

\begin{prop}\label{prop:existunique}
{\rm (1)} For any subspace $W\subset\mathbb{R}^N$ of dimension $N-d$, there exists a $d$-dimensional tight frame $\mathcal{S}$ of size $N$ such that $W=W(\mathcal{S})$.

{\rm (2)}~ 
Two tight frames $\mathcal{S}_1, \mathcal{S}_2$ are congruent if and only if $W(\mathcal{S}_1)=W(\mathcal{S}_2)$.

\end{prop}
\noindent{\it Proof}. (1)~ Let $p:\mathbb{R}^N\longrightarrow \mathbb{R}^N$ be the orthogonal projection onto $W^{\perp}$. Choosing an isometry $i:W^{\perp}\longrightarrow \mathbb{R}^d$, we put $P=ip$, which is obviously a frame projection satisfying (3) in Proposition \ref{prop:nnn} such that $W={\rm ker}~\!P$.
\smallskip 

(2)~ 
If ${\rm ker}~\!P_{\mathcal{S}_1}= {\rm ker}~\!P_{\mathcal{S}_2}(=W)$, then since $P_i=P_{\mathcal{S}_i}|W^{\perp}:W^{\perp}\longrightarrow\mathbb{R}^d$ is isometry, $U=P_1P_2{}^{-1}$ is an orthogonal transformation such that $P_{\mathcal{S}_1}=UP_{\mathcal{S}_2}$.
\hfill $\Box$

\medskip

The symmetric group $\mathfrak{S}_N$ of $\{1,2,\ldots,N\}$ acts on $\mathbb{R}^N$ as axis permutations\footnote{This is called the {\it standard representation} of $\mathfrak{S}_N$.}, i.e., 
$$
\sigma({\bf f}_i)={\bf f}_{\sigma(i)}\quad (\sigma\in \mathfrak{S}_N).
$$
The {\it automorphism group} ${\rm Aut}(\mathcal{S})$ of a 1-tight frame $\mathcal{S}=\{{\bf v}_i\}_{i=1}^N$ is defined to be the subgroup of $\mathfrak{S}_N$ consisting of $\sigma\in \mathfrak{S}_N$ satisfying $\sigma\big(W(\mathcal{S})\big)=W(\mathcal{S})$. By virtue of Proposition \ref{prop:existunique} (2), there is an injective homomorphism $U:{\rm Aut}(\mathcal{S})\longrightarrow O(d)$ such that $U(\sigma){\bf v}_i={\bf v}_{\sigma(i)}$. 

If ${\rm Aut}(\mathcal{S})$ acts transitively on $\{1,2,\ldots,N\}$, then $\mathcal{S}$ is said to be {\it isotropic}\footnote{An isotropic tight frame $\{{\bf v}_i\}_{i=1}^N$ is {\it uniform} in the sense that $\|{\bf v}_1\|=\cdots=\|{\bf v}_N\|$. The notion of uniform tight frame appears in various applications. The classification of isotropic frames is obviously related to that of subgroups of $\mathfrak{S}_N$ acting transitively on  $\{1,2,\ldots,N\}$ which has been pursued for over a century since 
the 1860 Grand Prix of the Acad\'{e}mie des Sciences.}. 

Tight frames associated with regular polygons, Platonic solids, and Archimedean solids are isotropic. Among them, the equilateral triangle and the regular tetrahedron are very special in the sense that ${\rm Aut}(\mathcal{S})$ agrees with $\mathfrak{S}_N$. We shall say that $\mathcal{S}$ with the property ${\rm Aut}(\mathcal{S})=\mathfrak{S}_N$ is {\it strongly isotropic}.

The higher dimensional analogue of the equilateral triangle and the regular tetrahedron is the {\it equilateral simplex}, a simplex $\Delta$ in $\mathbb{R}^d$ 
whose edges have equal length. Suppose that its barycenter is the origin $O$, and let $P_1,\ldots,P_{d+1}$ be its vertices. The symmetric group $\mathfrak{S}_{d+1}$ acts on $\Delta$ as orthogonal transformations. Thus $\{{\bf v}_i\}_{i=1}^{d+1}$ defined by ${\bf v}_i=\overrightarrow{OP_i}$ is a strongly isotropic tight frame. 

Conversely we have the following.

\begin{prop}\label{prop:equi}
Let $\mathcal{S}$ be a strongly isotropic 1-tight frame in $\mathbb{R}^d$ ($d\geq 2$). Then $\mathcal{S}$ is the frame associated with the equilateral simplex.
\end{prop}

\noindent{\it Proof}.~Since ${\rm Aut}(\mathcal{S})=\mathfrak{S}_N$, invariant subspaces for the ${\rm Aut}(\mathcal{S})$-action on $\mathbb{R}^N$ are either $W=\{{}^t(x,x,\ldots,x)|x\in \mathbb{R}\}$ or $W^{\perp}$ (indeed, $\mathbb{R}^N=W\oplus W^{\perp}$ gives the irreducible decomposition for the $\mathfrak{S}_N$-action). Because $d\geq 2$, the vanishing subspace for $\mathcal{S}$ must be $W$, and $N=d+1$. Then the frame projection $P$ is identified with the orthogonal projection of $\mathbb{R}^{d+1}$ onto $W^{\perp}$. Hence
$$
{\bf v}_i=P({\bf f}_i)={\bf f}_i-\frac{1}{d+1}\sum_{i=1}^{d+1}{\bf f}_i.
$$
Since $
\langle {\bf v}_i,{\bf v}_j\rangle=\delta_{ij}-\frac{1}{d+1}
$, if  ${\bf v}_i=\overrightarrow{OP_i}$, then $
(P_iP_j)^2=\|{\bf v}_i-{\bf v}_j\|^2=2 \quad (i\neq j).
$ 
  We thus conclude that $P_1,\ldots,P_{d+1}$ is the vertices of the equilateral simplex.
\hfill $\Box$

\medskip
The proof of the following proposition is left as an exercise for the reader.

\begin{prop}
If  $\{\overrightarrow{OP_i}\}_{i=1}^{d+1}$ is a tight frame of $\mathbb{R}^d$, and $\displaystyle\sum_{i=1}^{d+1}\overrightarrow{OP_i}={\bf 0}$, then $P_1,P_2,\ldots, P_{d+1}$ be the vertices of an equilateral simplex.
\end{prop}

We go back to the general case. 
In terms of matrices, what we have said in Proposition \ref{prop:nnn} is rephrased as 

\begin{prop}\label{prop:stief}
The row vectors in a matrix
$$ A=
\begin{pmatrix}
a_{11} & a_{12} & \cdots & a_{1d}\\
a_{21} & a_{22} & \cdots & a_{2d}\\
       &        & \cdots &       \\
       &        & \cdots &       \\
a_{N1} & a_{N2} & \cdots & a_{Nd}
\end{pmatrix}\in M_{N,d}(\mathbb{R})
$$
give rise to a 1-tight frame if and only if the column vectors of $A$ form an orthonormal system (i.e. ${}^t\!AA=I_d$).
\end{prop}

For if we define $\mathcal{S}=\{{\bf v}_i\}_{i=1}^N$ by writing ${}^t\!A=({\bf v}_1,\ldots,{\bf v}_N)$, then ${}^t\!A$ is the matrix for $P_{\mathcal{S}}$. 
 
Writing $A=({\bf a}_1,\ldots,{\bf a}_d)$, we find that ${\rm ker}~\!P=\langle{\bf a}_1,\ldots,{\bf a}_d\rangle_{\mathbb{R}}{}^{\perp}$ because 
$$
\langle{\bf a}_1,\ldots,{\bf a}_d\rangle_{\mathbb{R}}{}^{\perp}=({\rm Image}~\!A)^{\perp}={\rm ker}~\!{}^t\!A.
$$ 
Thus we have the following proposition which rephrases Proposition \ref{prop:existunique}.

\begin{prop}\label{prop:matrixrep}  Given an $(N-d)$-dimensional subspace $W$ of $\mathbb{R}^N$, there exists a solution $A\in M_{N,d}(\mathbb{R})$ of the equations
\begin{eqnarray}\label{eq:quadraticmatrix}
&&{}^t\!AA=I_d,\\
&& {}^t\!A{\bf x}={\bf 0} \quad({\bf x}\in W). 
\end{eqnarray}
If $A_1,A_2$ are solutions, then there exists $U\in O(d)$ such that $A_1=A_2U$.
\end{prop}

We now give an explicit parameterization of congruence classes of 1-tight frames. We denote by ${\rm T}_N(\mathbb{R}^d)$ the set of congruence classes of 1-tight frames $\{{\bf v}_i\}_{i=1}^N$ in $\mathbb{R}^d$. Proposition \ref{prop:stief} tells us that the set of $1$-tight frames is identified with the {\it Stiefel manifold} ${\rm V}_d(\mathbb{R}^N)$ $(=\{A\in M_{N,d}(\mathbb{R})|~{}^t\!AA=I_d\}$). The action of the orthogonal group $O(d)$ on ${\rm V}_d(\mathbb{R}^N)$ by $A\in {\rm V}_d(\mathbb{R}^N) \mapsto AU^{-1}\in {\rm V}_d(\mathbb{R}^N)$ \big($U\in O(d)$\big) is compatible with the action of $O(d)$ on the set of 1-tight frames because 
$$
{}^t(AU^{-1})=U~\!{}^t\!A=U({\bf v}_1,\ldots,{\bf v}_N)=\big(U({\bf v}_1),\ldots,U({\bf v}_N)\big).
$$
Therefore ${\rm T}_N(\mathbb{R}^d)$ is identified with the quotient space ${\rm V}_d(\mathbb{R}^N)/O(d)$, where the canonical projection $\varphi: {\rm V}_d(\mathbb{R}^N)\longrightarrow {\rm V}_d(\mathbb{R}^N)/O(d)={\rm T}_N(\mathbb{R}^d)$ coincides with the map which brings 1-tight frames to their congruence classes.

The quotient space ${\rm V}_d(\mathbb{R}^N)/O(d)$ is nothing but the the {\it Grassmannian} ${\rm Gr}_{N-d}($ $\mathbb{R}^N)$, i.e., the set of $(N-d)$-dimensional subspaces of $\mathbb{R}^N$. Therefore ${\rm T}_N(\mathbb{R}^d)$ is identified with ${\rm Gr}_{N-d}(\mathbb{R}^N)$, which is also identified with ${\rm Gr}_d(\mathbb{R}^N)$ via the correspondence $W\mapsto W^{\perp}$.
Under these identifications, the canonical projection $\varphi$ turns out to be nothing but the map ${\rm V}_d(\mathbb{R}^N)\longrightarrow {\rm Gr}_d(\mathbb{R}^N)$ giving the well-known structure of an $O(d)$-principal bundle.

If we ignore the order of vectors in tight frames, it is natural to take up the quotient space $\mathfrak{S}_N\backslash {\rm T}_N(\mathbb{R}^d)$ $\big(=\mathfrak{S}_N\backslash {\rm V}_d\big(\mathbb{R}^N)/O(d)\big)$ where
the symmetric group $\mathfrak{S}_N$ acts on ${\rm V}_d(\mathbb{R}^N)$ by $(a_{ij})\mapsto (a_{\sigma^{-1}(i)j})$ $(\sigma \in \mathfrak{S}_N)$. For a 1-tight frame $\mathcal{S}=\{{\bf v}_i\}_{i=1}^N$, the isotropy group of the point ${}^t({\bf v}_1,\ldots,{\bf v}_N)O(d)\in {\rm V}_d(\mathbb{R}^N)/O(d)$ coincides with ${\rm Aut}(\mathcal{S})$.

In view of Proposition \ref{prop:matrixrep}, the Stiefel manifold ${\rm V}_d(\mathbb{R}^N)$ is regarded as a ``quadric" in $M_{N,d}(\mathbb{R})$, and a 1-tight frame with the vanishing group $W$ is obtained as a point in the intersection of the quadric ${\rm V}_d(\mathbb{R}^N)$ and the subspace $\big\{A\in M_{N,d}(\mathbb{R})|~{}^tA{\bf x}={\bf 0}~({\bf x}\in W)\big\}$.
Such locution turns out to become more natural when we consider the set of ``oriented" congruence classes of 2-dimensional tight frames (see Sect.~\!\ref{sec:Parameterizations of crystallographic tight frames}). Here the set of orientated congruence classes of 1-tight frames is $\widetilde{{\rm T}}_N({\mathbb{R}^d})={\rm V}_d(\mathbb{R}^N)/SO(d)$, which is identified with the {\it oriented Grassmannian} $\widetilde{\rm Gr}_d(\mathbb{R}^N)$,  the manifold consisting of all oriented $d$-dimensional subspaces of $\mathbb{R}^N$. It is a double cover over ${\rm Gr}_d(\mathbb{R}^N)$.

\medskip

We close this section by giving a simple remark on Gramm matrices associated with frames.
 
\begin{prop}
Two frames $\mathcal{S}_1$ and $\mathcal{S}_2$ are congruent if and only if $G_{\mathcal{S}_1}=G_{\mathcal{S}_2}$. 
\end{prop}

\noindent {\it Proof}.~ It suffices to show that if ${\bf b}_1,\ldots,{\bf b}_N$ span $\mathbb{R}^d$, and 
$$
\langle{\bf b}_i,{\bf b}_j\rangle=\langle{\bf c}_i,{\bf c}_j\rangle\quad (i,j=1,\ldots, N)
$$ 
for vectors ${\bf c}_1,\ldots,{\bf c}_N$ in $\mathbb{R}^d$, 
then there exists an orthogonal matrix $U\in O(d)$ such that ${\bf c}_i=U{\bf b}_i$. Certainly this is true when ${\bf b}_1,\ldots,{\bf b}_N$ ($N=d$) is a basis.
In the general case, we take a basis ${\bf b}_{i_1},\ldots,{\bf b}_{i_d}$, and $U\in O(d)$ with ${\bf c}_{i_k}=U{\bf b}_{i_k}$. Then
$$
\langle {\bf c}_{i_h}, {\bf c}_k\rangle=\langle {\bf b}_{i_h}, {\bf b}_k\rangle
=\langle U{\bf b}_{i_h},U{\bf b}_k\rangle=\langle {\bf c}_{i_h}, U{\bf b}_k\rangle,
$$
from which it follows that ${\bf c}_k=U{\bf b}_k$ for every $k$, as required. \hfill $\Box$

\medskip
Proposition \ref{prop:nnn} (4) combined with the above proposition tells us that ${\rm T}_N(\mathbb{R}^d)$ is identified with 
$$
\{G\in M_{N}(\mathbb{R})|~{}^tG=G,~ G^2=G,~ {\rm rnak}~\!G=d\}.
$$

\section{Crystallographic tight frames}
We come now to the proper subject of this paper. 
Recall that the tight frame $\{\mathbf{v}_1,\mathbf{v}_2,\mathbf{v}_3\}$ mentioned in Introduction forms a crystalline pattern if and only if it generates a lattice in $\mathbb{R}^2$. Having this fact in mind, we shall introduce the notion of {\it crystallographic tight frame}. 

Before going into the subject, we review some items in the theory of lattice groups which are often used in the rest of this article.
In general, for a lattice $\mathcal{L}$ in an $n$-dimensional vector space $M$ with an inner prodcut $\langle\cdot,\cdot\rangle$, we denote by ${\rm vol}(M/\mathcal{L})$ the volume of the flat torus $M/\mathcal{L}$. If $\{{\bf a}_1,\ldots,{\bf a}_n\}$ is a $\mathbb{Z}$-basis of $\mathcal{L}$, then ${\rm vol}(M/\mathcal{L})^2=\det (\langle {\bf a}_i,{\bf a}_j\rangle)$. Further we have 
$$
{\rm vol}(M/\mathcal{L}_1)=|\mathcal{L}_2/\mathcal{L}_1|{\rm vol}(M/\mathcal{L}_2)
$$
for two lattices $\mathcal{L}_1, \mathcal{L}_2$ with $\mathcal{L}_1\subset\mathcal{L}_2$.

We denote by $\mathcal{L}^{\#}$ the {\it dual lattice} of $\mathcal{L}$; namely $\mathcal{L}^{\#}=\big\{{\bf x}\in M|~\langle {\bf x},{\bf y}\rangle\in \mathbb{Z}~({\bf y}\in \mathcal{L})\big\}$. We observe
$$
{\rm vol}(M/\mathcal{L}^{\#})={\rm vol}(M/\mathcal{L})^{-1}.
$$ 
Therefore if $\mathcal{L}$ is {\it integral}, i.e., $\mathcal{L}\subset\mathcal{L}^{\#}$, then 
\begin{equation}\label{eq:dualint}
|\mathcal{L}^{\#}/\mathcal{L}|={\rm vol}(M/\mathcal{L})^2,
\end{equation}
in particular, ${\rm vol}(M/\mathcal{L})^2$ is an integer.

\medskip



A frame $\mathcal{S}=\{\mathbf{v}_i\}_{i=1}^N$ of $\mathbb{R}^d$ is said to be {\it crystallographic} if it generates a lattice in $\mathbb{R}^d$, or what is the same is that $\mathcal{L}_{\mathcal{S}}={\rm Image}~\!\rho_{\mathcal{S}}$ is a lattice, where $
\rho_{\mathcal{S}}=P_{\mathcal{S}}|\mathbb{Z}^N:\mathbb{Z}^N\longrightarrow \mathbb{R}^d$. We shall designate $\mathcal{L}_{\mathcal{S}}$ as the {\it periodic lattice} for $\mathcal{S}$.

The {\it vanishing subgroup} $H(\mathcal{S})$ associated with a frame $\mathcal{S}$ is defined to be the kernel of the homomorphism 
$
\rho_{\mathcal{S}}
$, i.e., ${}^t(k_1,\ldots,k_N)\in H(\mathcal{S})$ if and only if
$
k_1\mathbf{v}_1+\cdots+k_N\mathbf{v}_N={\bf 0}
$. 

We have a sufficient condition for a frame being crystallographic (this turns out to be a necessary condition if the frame is 1-tight; see Proposition \ref{prop:necessary}).
 
\begin{lemma}\label{lemma:lattice}
A frame $\mathcal{S}=\{\mathbf{v}_i\}_{i=1}^N$ is crystallographic provided that the Gramm matrix $G_{\mathcal{S}}$ is essentially rational, i.e., $\langle {\bf v}_i,{\bf v}_j\rangle\in \lambda\mathbb{Q}$ ($i,j=1,\ldots,N$) with $\lambda>0$.
\end{lemma}

\noindent{\it Proof}.~ Put $v_{ij}=\langle {\bf v}_i,{\bf v}_j\rangle$. Suppose that $\{\mathbf{v}_i\}_{i=1}^N$ does not generate a lattice. Then ${\bf 0}$ is an accumulation point in $\mathbb{R}^d$ of the subgroup 
$$
\displaystyle\Big\{\sum_{i=1}^Nk_i{\bf v}_i|~k_i\in \mathbb{Z}~ (i=1,\ldots,N)\Big\},
$$
so that one can find a sequence $\{x_{ni}\}_{n=1}^{\infty}\in \mathbb{Z}$ such that 
$$\displaystyle\sum_{i=1}^Nx_{ni}{\bf v}_i\neq {\bf 0}~ \text{and}~ \displaystyle\lim_{n\to \infty}\sum_{i=1}^Nx_{ni}{\bf v}_i={\bf 0}.
$$ 
On the other hand, 
$$
\Big\|\sum_{i=1}^Nx_{ni}{\bf v}_i\Big\|^2=\sum_{i,j=1}^Nv_{ij}x_{ni}x_{nj}.
$$
Take a positive integer $N$ such that $\lambda^{-1}Nv_{ij}\in \mathbb{Z}$.
Since $\displaystyle\Big\{\lambda^{-1}N\sum_{i,j=1}^Nv_{ij}x_ix_j|~x_i\in \mathbb{Z}\Big\}$ is discrete, we have a contradiction. \hfill $\Box$

\medskip

The period lattice $\mathcal{L}_{\mathcal{S}}$ is just the set of vertices of the network obtained as the projected image by the frame projection of the {\it hypercubic lattice}. Here hypercubic lattice means the net associated with the $N$-dimensional standard lattice $\mathbb{Z}^N$, a generalization of the square and cubic lattices.

The hypercubic lattice as a graph is the {\it Cayley graph} $X$ associated with the free abelian group $\mathbb{Z}^N$ with the set of generators $\{{\bf f}_1,\ldots,{\bf f}_N\}$ (remember that $\{{\bf f}_i\}$ is the fundamental basis of $\mathbb{R}^N$). Thus the quotient graph $X/H(\mathcal{S})$ by the natural $H(\mathcal{S})$-action on $X$ is the Cayley graph associated with the factor group $\mathbb{Z}^N/H(\mathcal{S})$ with the set of generators $\big\{{\bf f}_i+H(\mathcal{S})\big\}_{i=1}^N$. Furthermore $X$ is the maximal abelian covering graph over the $N$-{\it bouquet graph}, the graph with a single vertex and $N$ loop edges. The projected image of the hypercubic lattice can be thought of as a (periodic) realization of the abstract graph $X/H(\mathcal{S})$ \big(possibly having degenerate edges, multiple edges, and/or colliding vertices when realized in $\mathbb{R}^d$, like the square lattice in Fig.~\!\ref{fig:projection} (a)\big). The map (graph morphism) $\omega$ of $X$ onto $X/H(\mathcal{S})$ associated with the canonical homomorphism $\mathbb{Z}^N\longrightarrow \mathbb{Z}^N/H(\mathcal{S})$ turns out to be a {\it covering map}, and is compatible with the frame projection $P$ of the hypercubic lattice onto the projected image.

\begin{equation*} \begin{CD}
X @>\text{}>>\mathbb{R}^N\\
@V\text{$\omega$}VV @VV\text{$P$}V \\
X/H(\mathcal{S}) @>>\text{} > \mathbb{R}^d
\end{CD}
\end{equation*}

This story is much generalized in Sect.~\!\ref{sect:Standard crystal models} in terms of topological crystals and their standard realizations.

\medskip

There are plenty of sources of crystallographic tight frames.
\smallskip

(1)~ The (isotropic) tight frame associated with the $N$-regular polygon is crystallographic if and only if $N=3,4,6$. The 2D crystal pattern for $N=3,6$ is the regular triangular lattice. For $N=4$, we have the square lattice.

\smallskip

(2)~ Platonic solids which yield crystallographic tight frames are the tetrahedron, cube, and octahedron\footnote{The reason why restricted polygons and polyhedra appear in (1), (2) is derived from the following general fact: If a finite subgroup $G$ of $GL_d(\mathbb{Z})$ contains an element with order $n$, then $\varphi(n)\leq d$, where $\varphi(n)$ is the Euler function, the number of positive integers $k$ with $1\leq k\leq n$ and ${\rm gcd}(n,k)=1$. Thus for $d=2,3$, the possible order is $2,3,4$, or $6$.}. The tetrahedron and cube yield the net associated with the body-centered lattice \big(Fig.~\!\ref{fig:3d} (e)\big), while the crystal net corresponding to the octahedron is the cubic lattice.

\begin{figure}[htbp]
\begin{center}
\includegraphics[width=.6\linewidth]
{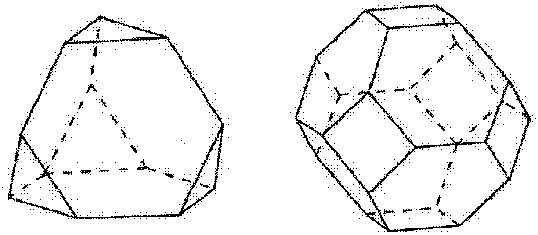}
\end{center}
\caption{Truncated tetrahedron and Truncated octahedron}\label{fig:archi}
\end{figure}

\begin{figure}[htbp]
\begin{center}
\includegraphics[width=.65\linewidth]
{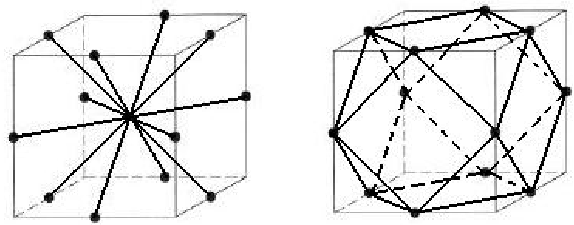}
\end{center}
\caption{$A_3$ and Cuboctahedron}\label{fig:cuboocta}
\end{figure}

Among all Archimedean solids, truncated tetrahedron, cuboctahedron, and truncated octahedron yield crystallographic tight frames; others do not\footnote{We may use Lemma \ref{lemma:lattice} to check this (see also Proposition \ref{prop:necessary}).} (Fig.~\!\ref{fig:archi} and \ref{fig:cuboocta}).

\smallskip

(3)~ 
An advanced example of tight frames is derived from root systems.
For the convenience of the reader, let us recall the definition \cite{hum}.  

A {\it root system} in $\mathbb{R}^d$ is a finite set $\varPhi$ of non-zero
vectors (called roots) that satisfy the following
conditions:

\smallskip

1. The roots span $\mathbb{R}^d$.

2. The only scalar multiples of a root $\mathbf{x}\in \varPhi$ that belong to $\varPhi$ are $\mathbf{x}$ itself and $-\mathbf{x}$.

3. For every root $\mathbf{x}\in \varPhi$, the set $\varPhi$ is closed under reflection through the hyperplane perpendicular to $\mathbf{x}$.

4. (Integrality) If $\mathbf{x}$ and $\mathbf{y}$ are roots in $\varPhi$, then the projection of $\mathbf{y}$ onto the line through $\mathbf{x}$ is a half-integral multiple
of $\mathbf{x}$.

\smallskip
The reflection $\sigma_{\bf x}$ through the hyperplane perpendicular to $\mathbf{x}$ is explicitly expressed as
$$
\sigma_{\bf x}({\bf y})={\bf y}-2\frac{\langle{\bf y},{\bf x} \rangle}{\|{\bf x}\|^2}{\bf x}.
$$

The group of orthogonal transformations of $\mathbb{R}^d$ generated by reflections through hyperplanes associated to the roots of $\varPhi$ is finite and called the {\it Weyl group} of $\varPhi$.

A root system $\varPhi$ is called {\it irreducible} if it cannot be partitioned into the union of proper subsets such that each root in one set is orthogonal to each root on the other. If $\varPhi$ is irreducible, then the Weyl group acts irreducibly on $\mathbb{R}^d$. Therefore in view of Proposition \ref{thm:tightirre}, $\varPhi$ (under any ordering of roots) gives a tight frame. $\varPhi$ is redundant in the sense that it contains both ${\bf x}$ and $-{\bf x}$. One can take subset $\varPhi^{+}$ such that $\varPhi=\varPhi^{+}\cup -\varPhi^{+}$ and $\varPhi^{+}\cap -\varPhi^{+}=\emptyset$ (the set of {\it positive roots} in the root system gives such $\varPhi^{+}$). When $\varPhi$ is irreducible, $\varPhi^{+}$ obviously gives a tight frame. 

The tight frames given by $\varPhi$ and $\varPhi^+$ are crystallographic because $\varPhi$ generates a lattice (called the {\it root lattice})\footnote{$\varPhi$ and $\varPhi^+$ yield the same crystal net.}.

It is known that there are four
infinite families of {\it classical} irreducible root systems designated as $A_d$ ($d\geq 1$), $B_d$ ($d\geq 2$), $C_d$ ($d\geq 3$), and $D_d$ ($d\geq 4$), and the five {\it exceptional} 
root systems $E_6$, $E_7$, $E_8$, $F_4$, and $G_2$. Among them, $A_d, D_d, E_6, E_7, E_8$ give isotropic frames (cf. \!\cite{Ebe}). The crystal nets for $A_2$ and $A_3$ are the regular triangular lattice and the net associated with the face-centered lattice \big(Fig.~\!\ref{fig:3d} (d)\big), respectively. The tight frame associated with $A_3$ coincides with the one coming from {\it cuboctahedron}  (Fig.~\!\ref{fig:cuboocta}). 

Interestingly $A_3$ (and the net associated with $A_3$) comes up in crystallography in various forms; say, in the hexagonal arrangement that gives the densest sphere packing\footnote{This is the Kepler conjecture for which Thomas Hales gave a proof in 1998. Johaness Kepler stated this conjecture in the short pamphlet entitled {\it  New-Year's gift concerning six-cornered snow} (``Strena Seu de Nive Sexangula" in Latin) in 1611.}, the diamond crystal, and the diamond twin (see Sect.~\!\ref{sect:finite graphs}).   

\begin{figure}[htbp]
\begin{center}
\includegraphics[width=.6\linewidth]
{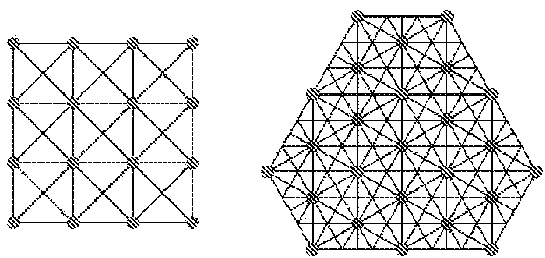}
\end{center}
\caption{The net for $B_2$ and $G_2$}\label{fig:B2G2}
\end{figure}

Figure \ref{fig:B2G2} depicts the nets for $B_2$ and $G_2$ (compare with Fig.~\!\ref{fig:coxeter}). 

(4)~ The {\it Leech lattice} $\Lambda$ discovered by John Leech in 1967 is a highly symmetrical lattice in $\mathbb{R}^{24}$ characterized by the following properties \cite{con}, \cite{milnor}:

$\cdot$~ It is self-dual, i.e., $\Lambda^{\#}=\Lambda$.

$\cdot$~ It is even, i.e., the square of the length of any vector in $\Lambda$ is an even integer.

$\cdot$~ The length of any non-zero vector in $\Lambda$ is at least 2.

\smallskip

The group (called the {\it Conway group} ${\rm Co}_0$) of orthogonal transformations which preserve $\Lambda$ permutes transitively the 196560 vectors ${\bf x}\in \Lambda$ with $\|{\bf x}\|=2$. The action of ${\rm Co}_0$ on $\mathbb{R}^{24}$ is irreducible, so that $\varPhi=\{{\bf x}\in \Lambda|~\|{\bf x}\|=2\}$ gives a crystallographic, isotropic tight frames. It is known that 
${\rm Co}_0$ modulo its center (designated as ${\rm Co}_1$) is a simple group of order  $4157776806543360000=2^{21}\cdot 3^9\cdot 5^4\cdot 7^2\cdot 11\cdot 13\cdot 23$. The Leech lattice is not a root lattice, but considered a kin of the exceptional root lattice $E_8$ in view of the fact that $E_8$ is the unique self-dual even lattice in 8-dimension.
\smallskip

\begin{figure}[htbp]
\begin{center}
\includegraphics[width=.45\linewidth]
{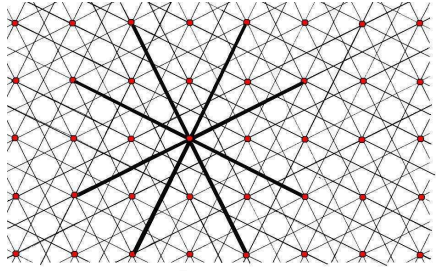}
\end{center}
\caption{A Pythagorean lattice}\label{fig:phyta}
\end{figure}

(5)~ Let $(x,y,z)$ be a primitive Pythagorean triple; namely $x,y,z$ are coprime positive integers satisfying $x^2+y^2=z^2$. Put
$$
\mathbf{v}_1=\begin{pmatrix}
1\\
0
\end{pmatrix},~
 \mathbf{v}_2=\begin{pmatrix}
x/z\\
y/z
\end{pmatrix},~ \mathbf{v}_3=\begin{pmatrix}
0\\
1
\end{pmatrix},~ \mathbf{v}_4=\begin{pmatrix}
-y/z\\
x/z
\end{pmatrix}.
$$
One can check that $\{\mathbf{v}_i\}_{i=1}^4$ is a crystallographic tight frame  whose vanishing subgroup is 
$\mathbb{Z}{}^t(z,-x,0,y)+\mathbb{Z}{}^t(0,y,-z,x)$. The crystalline pattern associated with this tight frame is what we call a {\it Pythagorean lattice} \cite{sunadadiamnod}. Note that the tight frame $\{\pm\mathbf{v}_i\}_{i=1}^4$ consisting 8 vectors is isotropic. See Fig.~\!\ref{fig:phyta} in the case $(x,y,z)=(3,4,5)$.

Related to Pythagorean lattices is the notion of {\it coincidence symmetry group}, which originates in the theory of crystalline interfaces and grain boundaries in polycrystalline materials. 

In general, the coincidence symmetry group  for a lattice $\mathcal{L}$ in $\mathbb{R}^d$ is defined as
$$
G(\mathcal{L})=\{g\in SO(d)|~g\mathcal{L}\sim \mathcal{L}\},
$$
where $\mathcal{L}_1\sim \mathcal{L}_2$ means that two lattices $\mathcal{L}_1$ and $\mathcal{L}_2$ are {\it commensurable}, i.e., $\mathcal{L}_1\cap \mathcal{L}_2$ is a lattice. Because $\sim$ is an equivalence relation, $G(\mathcal{L})$ is actually a subgroup of $SO(d)$. Note that the {\it symmetry group} $\{g\in SO(d)|~g\mathcal{L}= \mathcal{L}\}$ is always finite, while $G(\mathcal{L})$ could be infinite in general. In the special case $\mathcal{L}=\mathbb{Z}^d$, we get 
$G(\mathbb{Z}^d)=SO(d)\cap M_d(\mathbb{Q})$ which is a dense subgroup of $SO(d)$ (see the remark below). In particular 
$$
G({\bf Z}^2)=\Big\{\begin{pmatrix}
p & -q \\
q & p
\end{pmatrix}
\Big|~p^2+q^2=1, ~p,q\in \mathbb{Q}
\Big\}.
$$
If $p=x/z, q=y/z$ with a primitive Pythagorean triple $(x,y,z)$, then we have 
$$
g\mathbb{Z}^2=\mathbb{Z}\begin{pmatrix}
x/z\\
y/z\end{pmatrix}+\mathbb{Z}\begin{pmatrix}-y/z\\
x/z\end{pmatrix}.
$$

\medskip

(6)~ Let $\mathcal{L}^{\rm hc}$ be the 2-dimensional lattice with the ${\bf Z}$-basis ${\bf a}_1,{\bf a}_2$ such that $\|{\bf a}_1\|^2=\|{\bf a}_2\|^2=1$, ${\bf a}_1\cdot{\bf a}_2=-1/2$ (note that  $\mathcal{L}^{\rm hc}$ is a lattice whose translational action preserves the honeycomb). We then have
\begin{eqnarray*}
G(\mathcal{L}^{\rm hc})&=&\Big\{\begin{pmatrix}
p-\frac{1}{2}q & -\frac{\sqrt{3}}{2}q \\
\frac{\sqrt{3}}{2}q & p-\frac{1}{2}q
\end{pmatrix}
\Big|~p^2-pq+q^2=1, ~p,q\in \mathbb{Q}
\Big\}
\end{eqnarray*}
which is also a dense subgroup of $SO(2)$.

Put ${\bf v}_1={\bf a}_1$, ${\bf v}_2={\bf a}_2$, ${\bf v}_3=-{\bf a}_1-{\bf a}_2$. Obviously $\{{\bf v}_1,{\bf v}_2,{\bf v}_3\}$ is a 1-tight frame associated with the equilateral triangle whose period lattice is $\mathcal{L}^{\rm hc}$. For any $g\in G(\mathcal{L}^{\rm hc})$, the frame $\{{\bf v}_1,{\bf v}_2,{\bf v}_3, g{\bf v}_1,g{\bf v}_2,g{\bf v}_3\}$ is crystallographic and tight.

\medskip

Here are several remarks related to the examples (5), (6). 

\smallskip

(1)~ The coincidence symmetry group $G(\mathcal{L})$ is dense in $SO(d)$~ $\Leftrightarrow$~ the lattice $\mathcal{L}$ is {\it essentially rational}, i.e.,
there exists a positive $\lambda$ such that $\langle \mathcal{L},\mathcal{L}\rangle\subset \lambda\mathbb{Q}$. To sketch the proof, we first note that $G(\alpha\mathcal{L})=G(\mathcal{L})$ for any lattice and $\alpha>0$. Thus in the proof of the implication $\Leftarrow$, one may assume that $\mathcal{L}$ is rational. Selecting a $\mathbb{Z}$-basis $\{{\bf a}_1,\ldots,{\bf a}_d\}$ of $\mathcal{L}$, we consider the symmetric matrix $S={}^t({\bf a}_1,\ldots,{\bf a}_d)({\bf a}_1,\ldots,{\bf a}_d)$, and define 
$$
L_{\mathbb{Q}}(S)=\{X=(x_{ij})\in M_d(\mathbb{Q}),~{}^{t}\hspace{-0.05cm}XS+SX=O\},
$$
which is a vector space over $\mathbb{Q}$ of dimension $\frac{d(d-1)}{2}$ because $S$ is rational. 
Then putting
$$
\varphi(X)=({\bf a}_1,\ldots,{\bf a}_d)(I-X)(I+X)^{-1}({\bf a}_1,\ldots,{\bf a}_d)^{-1},
$$
we have an injective map $\varphi:L_{\mathbb{Q}}(S)\longrightarrow G(\mathcal{L})$ whose image is dense in $SO(d)$ ($\varphi$ is what we call {\it Cayley's parameterization}\footnote{The idea, originally given by A. Cayley, dates back to 1846.}). From this argument, it follows that $G(\mathcal{L})$ is dense in $SO(d)$. 

Conversely suppose that $G(\mathcal{L})$ is dense in $SO(d)$, which is equivalent to the condition that in the rotation group $SO(S)=\{A\in GL_d(\mathbb{R})|~{}^t\!ASA=S,~\det A>0\}$ for the symmetric matrix $S$,
the subgroup $SO_{\mathbb{Q}}(S)=\{A\in GL_d(\mathbb{Q})|~{}^t\!ASA=S,~\det A>0\}$ is dense. We easily find
$$
\big\{S'\in M_d(\mathbb{R})|~{}^t\!AS'A=S'~\big(A\in SO(S)\big),~{}^tS'=S'\big\}=\mathbb{R}S
$$
(this is equivalent to say that any $SO(d)$-invariant symmetric bilinear form on $\mathbb{R}^d$ is a scalar multiple of the standard inner product). On the other hand, the equation for $S'$ given by
$$
{}^t\!AS'A=S'~\big(A\in SO_{\mathbb{Q}}(S)\big),~{}^tS'=S'
$$
reduces to a homogeneous linear equation with rational coefficients having a non-zero real solution; say, $S$, so that one can find a non-zero rational solution $S_0$. Since $SO_{\mathbb{Q}}(S)$ is dense in $SO(S)$ by the assumption, we conclude that ${}^t\!AS_0A=S_0~\big(A\in SO(S)\big),~{}^tS_0=S_0$, thereby $S_0=\lambda S$ for some $\lambda\neq 0$, and $\mathcal{S}$ being essentially rational.

\smallskip

(2)~ For two frames $\mathcal{S}_1=\{{\bf u}_i\}_{i=1}^M$ and $\mathcal{S}_2=\{{\bf v}_i\}_{i=1}^N$, we define the {\it join} $\mathcal{S}_1 \vee \mathcal{S}_2=\{{\bf w}_i\}_{i=1}^{M+N}$ by 
$$
{\bf w}_i=\begin{cases}
{\bf u}_i & (i=1,\ldots,M)\\
{\bf v}_{i-M} & (i=M+1,\ldots,M+N).
\end{cases}
$$ 
If both $\mathcal{S}_1$ and $\mathcal{S}_2$ are tight frames, then so is $\mathcal{S}_1 \vee \mathcal{S}_2$. In order that $\mathcal{S}_1 \vee \mathcal{S}_2$ is crystallographic, it is necessary and sufficient that the period lattices $\mathcal{L}_{\mathcal{S}_1}$ and $\mathcal{L}_{\mathcal{S}_2}$ are commensurable. Obviously the period lattice of $\mathcal{S}_1 \vee \mathcal{S}_2$ is 
$\mathcal{L}_{\mathcal{S}_1}+\mathcal{L}_{\mathcal{S}_2}$.

\smallskip

(3)~ Since finite subgroups of $O(2)$ are cyclic or dihedral (cf.~\!\cite{new}), 2-dimensional crystallographic isotropic tight frames are classified as follows (note that the order of a cyclic subgroup must be less than or equal to $6$). 

\begin{figure}[htbp]
\begin{center}
\includegraphics[width=.65\linewidth]
{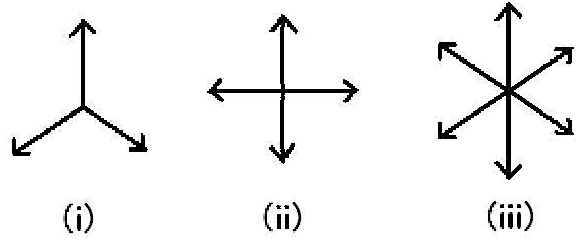}
\end{center}
\caption{2-dimensional crystallographic isotropic tight frames}\label{fig:iso}
\end{figure}

\smallskip
(a) Frames (i), (ii), (iii) depicted in Fig.~\!\ref{fig:iso}.

\smallskip

(b) Joins of each frame in (a) and its rotation by an element in the coincidence symmetry group \big(more precisely $G(\mathbb{Z}^2)$ for (ii) and $G(\mathcal{L}^{\rm hc})$ for (i) and (iii)\big) .

(4)~ In the terminology of crystallography, the intersection $\mathcal{L}_1\cap\mathcal{L}_2$ for commensurable lattices $\mathcal{L}_1$ and $\mathcal{L}_2$ is called the {\it coincidence site lattice} (CSL), while the sum $\mathcal{L}_1+\mathcal{L}_2$ is called the {\it displacement shift complete lattice} (DSC lattice). For more about the CSL theory, refer to \cite{grimmer}, \cite{baake}, \cite{zou}, \cite{zou1}.

\section{Parameterizations of crystallographic tight frames}\label{sec:Parameterizations of crystallographic tight frames}

If a frame $\mathcal{S}$ is crystallographic, then the vanishing group $H(\mathcal{S})$ is obviously a direct summand\footnote{A subgroup $B$ of an additive group $A$ is called a {\it direct summand} of $A$ if there exists another subgroup $B'$ such that $A=B\oplus B'$, i.e., any element $a\in A$ can be expressed uniquely as $a=b+b'$ with $b\in B, b'\in B'$. In the case $A=\mathbb{Z}^N$, $B$ is a direct summand if and only if the factor group $A/B$ is free abelian.} of $\mathbb{Z}^N$, and ${\rm rank}~\!H(\mathcal{S})=N-d$. Therefore $H(\mathcal{S})$ is a lattice in the vanishing subspace $W(\mathcal{S})={\rm ker}(P_{\mathcal{S}}:\mathbb{R}^N\rightarrow \mathbb{R}^d)$, and hence the subspace $H(\mathcal{S})_{\mathbb{R}}$ of $\mathbb{R}^N$ spanned by $H(\mathcal{S})$ coincides with $W(\mathcal{S})$. 

Conversely suppose that $H(\mathcal{S})$ is a direct summand of $\mathbb{Z}^N$.  Take a $\mathbb{Z}$-basis $\{{\bf c}_1,\ldots,$ ${\bf c}_N\}$ of $\mathbb{Z}^N$ such that $\{{\bf c}_{d+1},\ldots,{\bf c}_N\}$ is a $\mathbb{Z}$-basis of $H(\mathcal{S})$. Then $\big\{\rho_{\mathcal{S}}({\bf c}_1),$ $\ldots,\rho_{\mathcal{S}}({\bf c}_d)\big\}$ is a $\mathbb{Z}$-basis of $\mathcal{L}_{\mathcal{S}}$ (remember that $\rho_{\mathcal{S}}=P_{\mathcal{S}}|\mathbb{Z}^N$). Therefore $\mathcal{S}$ is crystallographic.

\begin{prop}\label{prop:Existence and Uniqueness}
{\rm (Existence and Uniqueness)}~ Let $H$ be a direct summand of $\mathbb{Z}^N$ of rank $N-d$. Then there exists a unique crystallographic 1-tight frame $\mathcal{S}$ of $\mathbb{R}^d$ up to congruence such that $H(\mathcal{S})=H$.

\end{prop}

\noindent {\it Proof}.\quad Let $H_{\mathbb{R}}$ be the subspace of $\mathbb{R}^N$ spanned by $H$, and let $P:\mathbb{R}^N\longrightarrow \mathbb{R}^d$ be the frame projection giving a 1-tight frame whose kernel is $H_{\mathbb{R}}$. 
To prove $H(\mathcal{S})=H$, note $H(\mathcal{S})={\rm ker}~\!P\cap\mathbb{Z}^N\supset H$ and ${\rm rank}~\!H(\mathcal{S})={\rm rank}~\!H$. Hence we have a finite subgroup $H(\mathcal{S})/H$ of $\mathbb{Z}^N/H$. Since $H$ is a direct summand, $H(\mathcal{S})/H=\{0\}$.  

The uniqueness follows from Proposition \ref{prop:existunique}. \hfill $\Box$

\medskip

There is a one-to-one correspondence between direct summands of rank $k$ in $\mathbb{Z}^N$ and subspaces of dimension $k$ in $\mathbb{Q}^N$. Indeed we have
\smallskip

{\rm (1)} Let $H\subset\mathbb{Z}^N$ be a direct summand, and let $H_{\mathbb{Q}}$ be the subspace of the vector space $\mathbb{Q}^N$ over $\mathbb{Q}$ spanned by $H$. Then $H_{\mathbb{Q}}\cap \mathbb{Z}^N=H$.

{\rm (2)} Let $W$ is a subspace of $\mathbb{Q}^N$ of dimension $k$. Then $H=W\cap \mathbb{Z}^N$ is a direct summand of $\mathbb{Z}^N$ of rank $k$.

\smallskip

In fact, (1) is proved in the same way as the proof of the above theorem. For (2), suppose that $\mathbb{Z}^N/(W\cap \mathbb{Z}^N)$ is not free. Then there exists $x\in \mathbb{Z}^N$ not contained in $W\cap \mathbb{Z}^N$ such that $nx\in W\cap \mathbb{Z}^N$. This implies that $x\in W$, and hence $x\in W\cap \mathbb{Z}^N$, thereby a contradiction.
\medskip

We denote by ${\rm T}_N^{\rm cr}(\mathbb{R}^d)$ the set of congruence classes of crystallographic 1-tight frames $\{{\bf v}_i\}_{i=1}^N$ in $\mathbb{R}^d$. Proposition \ref{prop:Existence and Uniqueness} together with the above remark tells us that ${\rm T}_N^{\rm cr}(\mathbb{R}^d)$ is parameterized by ${\rm Gr}_{N-d}(\mathbb{Q}^N)$. 


${\rm Gr}_{N-d}(\mathbb{Q}^N)$ is thought of as the set of ``rational points" on the complex Grassmannian ${\rm Gr}_{N-d}(\mathbb{C}^N)$ by considering it as a subvariety of the complex projective space $P^M(\mathbb{C})$ ($M=\binom{N}{d}-1$) by means of {\it Pl\"{u}cker embedding}. Here we recall the general definition of rational points. Let $K$ be an algebraic number field, and let $V$ be a projective algebraic variety, defined in some projective space $P^{n-1}(\mathbb{C})$ by homogeneous polynomials $f_1,\ldots,f_m$ with coefficients in $K$. A $K$-rational point of $V$ is a point $[z_1,\ldots,z_n]$ in $P^{n-1}(K)\big(\subset P^{n-1}(\mathbb{C})\big)$ that is a common solution of all the equations $f_j=0$.

Given $W\in {\rm Gr}_{N-d}(\mathbb{Q}^N)$, we shall explicitly construct a 1-tight frame whose vanishing group is $H=W\cap \mathbb{Z}^N$ (this is the stage where quadratic Diophantine equations show up). For this sake, we select non-zero vectors ${\bf n}_1,\ldots,{\bf n}_d$ of $W$ such that ${\bf n}_i\in \mathbb{Z}^N$, $\langle {\bf n}_i,{\bf x}\rangle=0$ for every ${\bf x}\in H$, and $\langle{\bf n}_i,{\bf n}_j\rangle=0$ for $i\neq j$ (we may do this because for a subspace $V$ of $\mathbb{Q}^N$, the ``rational" orthogonal complement $V^{\perp}=\big\{{\bf x}\in \mathbb{Q}^N|~\langle {\bf x},{\bf y}\rangle=0 ~({\bf y}\in V)\big\}$ satisfy $\mathbb{Q}^N=V\oplus V^{\perp}$; then we argue by induction for the construction of ${\bf n}_i$). Write $\|{\bf n}_i\|^2=m_i^2 D_i$ with a square free positive integer $D_i$ and a positive integer $m_i$. Then putting 
$
{\bf a}_i=(m_i\sqrt{D_i})^{-1}{\bf n}_i
$, 
we obtain $A=({\bf a}_1,\ldots,{\bf a}_d)\in {\rm V}_d(\mathbb{R}^N)$. If we define ${\bf v}_i$ by the relation $({\bf v}_1,\ldots,{\bf v}_N)={}^t\!A$, then we obtain a 1-tight frame $\mathcal{S}=\{{\bf v}_i\}_{i=1}^N$ in $\mathbb{R}^d$ whose vanishing group is $H$.

To express the matrix $A$ more explicitly, we write ${\bf n}_i={}^t(n_{1i},\ldots.n_{Ni})$. Then
\begin{equation}\label{eq:crystalpoint}
A=\begin{pmatrix}
\displaystyle\frac{1}{\sqrt{D_1}}\frac{n_{11}}{m_1}& \displaystyle\frac{1}{\sqrt{D_2}}\frac{n_{12}}{m_2} & \cdots &\displaystyle\frac{1}{\sqrt{D_d}}\frac{n_{1d}}{m_d}\\
\displaystyle\frac{1}{\sqrt{D_1}}\frac{n_{21}}{m_1}&\displaystyle \frac{1}{\sqrt{D_2}}\frac{n_{22}}{m_2} & \cdots&\displaystyle \frac{1}{\sqrt{D_d}}\frac{n_{2d}}{m_d}\\
 & \cdots & \cdots & \\
\displaystyle\frac{1}{\sqrt{D_1}}\frac{n_{N1}}{m_1}&\displaystyle \frac{1}{\sqrt{D_2}}\frac{n_{N2}}{m_2} & \cdots&\displaystyle \frac{1}{\sqrt{D_d}}\frac{n_{Nd}}{m_d}
\end{pmatrix},
\end{equation}
where $n_{ij}$ and $m_i$ satisfy the following Diophantine equation\footnote{A trivial remark is that a {\it ruler-compass construction} of the 1-tight frame is possible once this Diophantine equation would be solved.}
\begin{eqnarray*}
&&n_{1i}{}^2+\cdots+n_{Ni}{}^2=D_im_i{}^2 \quad (i=1,\ldots, d),\\
&&n_{1i}n_{1j}+\cdots+n_{Ni}n_{Nj}=0\quad (i\neq j).
\end{eqnarray*}  

Taking a look at (\ref{eq:crystalpoint}), we have

\begin{prop}\label{prop:necessary}
 The Gramm matrix for a crystallographic 1-tight frame is rational.
\end{prop}

Conversely if we start with a matrix $A\in {\rm V}_d(\mathbb{R}^N)$ of the form (\ref{eq:crystalpoint}), 
then the 1-tight frame associated with $A$ is crystallographic in view of Lemma \ref{lemma:lattice}. From this observation, it also follows that the set of congruence classes of 1-tight frames is parameterized by
$$
\{G\in M_N(\mathbb{Q})|~{}^tG=G,~ G^2=G,~ {\rm rank}~\!G=d \}.
$$

We also have

\begin{cor}\label{cor:rat}
The period lattice $\mathcal{L}_{\mathcal{S}}$ for a crystallographic 1-tight frame  $\mathcal{S}$ is rational in the sense that $\langle \mathcal{L}_{\mathcal{S}}, \mathcal{L}_{\mathcal{S}}\rangle \subset \mathbb{Q}$.

\end{cor}

We now restrict ourselves to the case $d=2$ and $N>2$, and consider $\widetilde{T}_N(\mathbb{R}^2)={\rm V}_2(\mathbb{R}^N)/SO(2)$, the set of oriented congruence classes of 1-tight frames in the plane $\mathbb{R}^2$. As remarked before, ${\rm V}_2(\mathbb{R}^N)/SO(2)$ is identified with the oriented Grassmannian $\widetilde{\rm Gr}_2(\mathbb{R}^N)$, which is also identified with the complex quadric
$$
Q_N=\{[z_1,\ldots,z_N]\in P^{N-1}(\mathbb{C})|~z_1{}^2+\cdots+z_N{}^2=0
\}.
$$
This routine procedure is carried out by the map
$$
\begin{pmatrix}
a_1 & b_1\\
a_2 & b_2\\
\cdot&\cdots\\
\cdot&\cdot\\
a_N & b_N
\end{pmatrix}\in {\rm V}_2(\mathbb{R}^N) \mapsto (a_1+b_1\sqrt{-1},\ldots,a_N+b_N\sqrt{-1}) \in \mathbb{C}^N,
$$
whichi 
yields a one-to-one correspondence between  $\widetilde{\rm T}_N(\mathbb{R}^2)$ and $Q_N$ \big(note that the quotient of $Q_N$ by the conjugation $[z_1,\ldots,z_N]\mapsto [\overline{z_1},\ldots,\overline{z_N}]$ is identified with ${\rm T}_N(\mathbb{R}^2)$\big).

For a direct summand $H$ of rank $N-d$ in $\mathbb{Z}^N$, define the subspace  $L_H$ in $P^{N-1}(\mathbb{C})$ by setting
\begin{eqnarray*}
L_H&=&\{[z_1,\ldots,z_N]\in P^{N-1}(\mathbb{C})|~k_1z_1+\cdots+k_Nz_N=0\\
&& \qquad \qquad \qquad \quad \text{for every}~{}^t(k_1,\ldots,k_N)\in H\}.
\end{eqnarray*}
The intersection $Q_N\cap L_H$ consists of two points (one is the conjugate of another) both of which correspond to an oriented similarity class of a crystallographic tight frame $\mathcal{S}$ with the vanishing subgroup $H=H_{\mathcal{S}}$. 

Applying the above observation in the general case, one can easily show that 
there exists a square free positive integer $D$ such that the two points in $Q_N\cap L_H$ are $\mathbb{Q}(\sqrt{-D})$-rational points on $Q_N$. 
We thus have

\begin{thm}\label{thm:2dimension}
The set 
$$
Q_N\cap \bigcup_{D}P^{N-1}\big(\mathbb{Q}(\sqrt{-D})\big)
$$ 
is identified with the set of oriented similarity classes of 2-dimensional tight frames of size $N$. 

\end{thm}

For an illustration, we shall make a brief excursion into $\mathbb{Q}(\sqrt{-D})$-rational points on $Q_3$. These rational points are related to 2D crystal patterns obtained as projected images of the cubic lattice mentioned in Introduction. 

The crystallographic 1-tight with the vanishing group $H=\mathbb{Z}~\!{}^t(n_1,n_2,n_3)$ corresponds to the point $[z_1,z_2,z_3]\in P^2({\bf C})$ 
given by
\begin{eqnarray*}
&& z_1=n_2{}^2+n_3{}^2,\\
&&z_2=-n_1n_2\pm \sqrt{-(n_1{}^2+n_2{}^2+n_3{}^2)},\\
&&z_3=-n_1n_3\mp \sqrt{-(n_1{}^2+n_2{}^2+n_3{}^2)}.
\end{eqnarray*}
Therefore $D$ is the square free part of $n_1{}^2+n_2{}^2+n_3{}^2$. 

A question arises: For which square free $D$, does the quadric $Q_3$ have a $\mathbb{Q}(\sqrt{-D})$-rational point? The answer is given in the following.

\begin{prop}\label{prop1} The quadric $Q_3$ has a $\mathbb{Q}(\sqrt{-D})$-rational point if and only if $D$ is not of the form $8k+7$. 

\end{prop}

This is, as easily conceived and proved below, a consequence of the theorem of three squares due to Legendre\footnote{Legendre' proof is not complete. It is Gauss who gave a complete proof.} (1798) which says that a positive integer $n$ can be expressed as the sum of three squares if and only if $n$ is not of the form $4^{\ell}(8k+7)$. 

From what we observed above, it follows that the quadric $Q_3$ has a $\mathbb{Q}(\sqrt{-D})$-rational point if and only if the equation $n_1{}^2+n_2{}^2+n_3{}^2=Dm^2$ has a non-trivial integral solution $n_1,n_2,n_3,m$.

We first show that if $D$ is not of the form $8k+7$, then $n_1{}^2+n_2{}^2+n_3{}^2=D$ has an integral solution $n_1,n_2,n_3$. To this end, suppose that $n_1{}^2+n_2{}^2+n_3{}^2=D$ has no integral solution. Then $D=4^{\ell}(8k+7)$ for some $\ell$ and $k$ by invoking Legendre's theorem. Since $D$ is square free, $D$ must be of the form $8k+7$. 

Next suppose that $D$ is of the form $8k+7$. If $n_1{}^2+n_2{}^2+n_3{}^2=Dm^2$ has a non-trivial integral solution $n_1,n_2,n_3,m$, then writing $m=2^{\ell}(2s+1)$, we have
$$
Dm^2=4^{\ell}(8k+7)(2s+1)^2=4^{\ell}\big(
8(8kt+7t+k)+7
\big),
$$  
where $t=s(s+1)/2$. This contradicts to Legendre's theorem.

\section{A height function on the rational Grassmannian}

The rational Grassmannian ${\rm Gr}_{N-d}(\mathbb{Q}^N)$ is dense in  ${\rm Gr}_{N-d}(\mathbb{R}^N)$; in particular there are infinitely many congruence classes of $d$-dimensional crystallographic 1-tight frames of size $N$. This section is devoted to a brief explanation how to count congruence classes, with some excursions into a characterization of crystallographic tight frames by means of a certain minimal principle. The tool that we employ is a natural height function on the rational Grassmannian.

We define the {\it height function} $h$ on ${\rm Gr}_k(\mathbb{Q}^N)$ in the following way. For $W\in {\rm Gr}_k(\mathbb{Q})$, we let $H=W\cap\mathbb{Z}^N$, which is a direct summand of $\mathbb{Z}^N$ as noticed before, and is an integral lattice in $H_{\mathbb{R}}(=W_{\mathbb{R}})$. We then put
$$
h(W)={\rm vol}(H_{\mathbb{R}}/H)\big(=|H^{\#}/H|^{1/2}\big),
$$
where the inner product on $H_{\mathbb{R}}$ is the one induced from the standard inner product on $\mathbb{R}^N$. 

The function $h$ deserves to be called a height function because we have

\begin{prop}\label{prop:heightfunction} For any $c>0$, there are only finitely many $W\in {\rm Gr}_k(\mathbb{Q}^N)$ such that $h(W)\leq c$.
\end{prop}

\noindent{\it Proof}.~ It suffices to prove that for every $c>0$, there are only finitely many direct summands $H$ of rank $k$ such that ${\rm vol}(H_{\mathbb{R}}/H)\leq c$.
To this end, we shall introduce the quantity $c(H)$. 
For a $\mathbb{Z}$-basis $J=\{{\bf c}_{1},\ldots, {\bf c}_k\}$ of $H$, we put
$$
c(J)=\max (\|{\bf c}_{1}\|,\ldots, \|{\bf c}_k\|),
$$ 
and define 
$$
c(H)=\min_{J}c(J).
$$
It is straightforward to check that there are only finitely many direct summands $H$ of rank $k$ such that $c(H)\leq c$. Therefore it suffices to prove 
$$
Cc(H)\leq {\rm vol}(H_{\mathbb{R}}/H)
$$
with a positive constant $C$ not depending on $H$. For this sake, we invoke the fact that there exists a $\mathbb{Z}$-basis (called a {\it reduced basis}) ${\bf u}_1,\ldots,{\bf u}_k$ of $H$ satisfying

\smallskip

(i)~ $\|{\bf u}_1\|\leq \|{\bf x}\|$ for all ${\bf x}\in H\backslash\{{\bf 0}\}$,

\smallskip

(ii) if ${\bf u}_1,\ldots,{\bf u}_{i-1},{\bf x}$ is a part of a $\mathbb{Z}$-basis of $H$, then $\|{\bf u}_i\|\leq \|{\bf x}\|$.

\smallskip
\noindent This property of ${\bf u}_1,\ldots,{\bf u}_k$ implies $\|{\bf u}_1\|\leq \|{\bf u}_2\|\leq \cdots
\leq \|{\bf u}_k\|$. Moreover a theorem 
 in the geometry of numbers due to Minkowski asserts that there exists a positive constant $C_k$ such that 
$$
\|{\bf u}_1\|\cdots\|{\bf u}_k\|\leq C_k{\rm vol}(H_{\mathbb{R}}/H)
$$
(cf.~\!\cite{borel}). By the definition of $c(J)$, if $J=\{{\bf u}_{1},\ldots, {\bf u}_k\}$, then $c(J)=\|{\bf u}_k\|$, so we find $c(H)\leq \|{\bf u}_k\|$. Since $\|{\bf u}_i\|\geq 1$ for $i=1,\ldots,k-1$ (note $H\subset \mathbb{Z}^N$), we have
$$
c(H)\leq \|{\bf u}_k\|\leq \|{\bf u}_1\|\cdots\|{\bf u}_k\|\leq C_k{\rm vol}(H_{\mathbb{R}}/H).
$$
This completes the proof. 
\hfill $\Box$

\begin{remark} {\rm (1)}~ By Hadamard' inequality, we have
$$
{\rm vol}(H_{\mathbb{R}}/H)
\leq \|{\bf c}_1\|\cdots \|{\bf c}_k\|\leq c(J)^k
$$
for any $\mathbb{Z}$-basis $J=\{{\bf c}_1,\ldots,{\bf c}_k\}$; therefore 
$$
{\rm vol}(H_{\mathbb{R}}/H)\leq c(H)^k.
$$

\smallskip

{\rm (2)}~ For a direct summand $H=\mathbb{Z}{}^t(n_1,\ldots,n_N)$ of $\mathbb{Z}^N$, 
$$
h(H_{\mathbb{Q}})=\sqrt{n_1{}^2+\cdots+n_N{}^2}.
$$ 
Using this, we obtain the asymptotic formula 
$$
\big|\{W\in P^{N-1}(\mathbb{Q})={\rm Gr}_1(\mathbb{Q}^N)~| h(W)\leq h\} \big|\sim \frac{1}{2}\zeta(N)^{-1}\omega_Nh^N\quad (h\to\infty),
$$ 
where $\omega_N=\pi^{N/2}/\Gamma(1+N/2)$, the volume of the unit disk in $\mathbb{R}^N$, and $\zeta(s)$ is the Riemann zeta function. It is interesting to ask whether a similar asymptotic formula holds for ${\rm Gr}_k(\mathbb{Q}^N)$. 

\end{remark}

The height function $h(W)$ is closely connected with crystallographic tight frames. 

\begin{prop}\label{prop:heightvolume} 
${\rm vol}(\mathbb{R}^d/\mathcal{L}_{\mathcal{S}})=h(H_{\mathbb{Q}})^{-1}\big(={\rm vol}(H_{\mathbb{R}}/H)^{-1}\big)$ for a crystallographic 1-tight frame $\mathcal{S}=\{{\bf v}_i\}_{i=1}^N$ whose vanishing subgroup is $H$.
\end{prop}

\noindent{\it Proof}.~  
Let $\mathcal{S}=\{{\bf v}_i\}_{i=1}^N$ be a crystallographic 1-tight frame of $\mathbb{R}^d$, and let $P=P_{\mathcal{S}}:\mathbb{R}^N\longrightarrow\mathbb{R}^n$ be the frame projection.
Again take a $\mathbb{Z}$-basis $\{{\bf c}_1,\ldots,$ ${\bf c}_N\}$ of $\mathbb{Z}^N$ as before such that $\{{\bf c}_{d+1},\ldots,{\bf c}_N\}$ is a $\mathbb{Z}$-basis of the vanishing subgroup $H(\mathcal{S})$. Note that $\big|\det({\bf c}_1,\ldots,{\bf c}_N)\big|=1$ because $({\bf c}_1,\ldots,{\bf c}_N)\in GL_{N}(\mathbb{Z})$. 
Thus for the square matrix $C=\big(\langle {\bf c}_i,{\bf c}_j\rangle\big)$ with integral entries, we find that 
$
\det C=\big(\det({\bf c}_1,\ldots,{\bf c}_N)\big)^2=1
$.

Putting ${\bf b}_i=P({\bf c}_i)$ ($i=1,\ldots,d$), we have a $\mathbb{Z}$-basis $\{{\bf b}_1,\ldots,{\bf b}_d\}$ of the period lattice $\mathcal{L}_{\mathcal{S}}$. Using ${}^t\!PP({\bf c}_i)-{\bf c}_i\in {\rm ker}~\!P$, 
we may write
$$
{}^t\!P({\bf b}_i)-{\bf c}_i=\sum_{j=d+1}^{N}f_{ij}{\bf c}_j.
$$ 
Substituting ${}^t\!PP({\bf c}_i)={\bf c_i}+\displaystyle\sum_{j=d+1}^{N}f_{ij}{\bf c}_j$ into $\langle {}^t\!PP({\bf c}_i),{\bf c}_k\rangle=0$ ($k=d+1,\ldots,b$), we have
$$
-\langle {\bf c}_i,{\bf c}_k\rangle=\sum_{j=d+1}^{N}f_{ij}\langle {\bf c}_j,{\bf c}_k \rangle.
$$
Now writing
$$
C=\begin{pmatrix}
C_{11} & C_{12}\\
C_{21} & C_{22}
\end{pmatrix},
$$
where $C_{11}\in M_d(\mathbb{R})$
and using these matrices of small size, we may compute the matrix $F=(f_{ij})\in M_{d,N-d}(\mathbb{R})$
as
$$
F=-C_{12}C_{22}^{-1}.
$$
Furthermore
\begin{eqnarray*}
\langle{\bf b}_i,{\bf b}_j\rangle&=&\langle{}^t\!PP({\bf c}_i),{\bf c}_j\rangle
=\Big\langle {\bf c}_i+\sum_{k=d+1}^Nf_{ik}{\bf c}_k,{\bf c}_j\Big\rangle\\
&=& \langle {\bf c}_i,{\bf c}_j\rangle+\sum_{k=d+1}^Nf_{ik}\langle{\bf c}_k,{\bf c}_j\rangle.
\end{eqnarray*}
Therefore the matrix $B=(\langle{\bf b}_i,{\bf b}_j\rangle)\in M_d(\mathbb{R})$  is computed as
\begin{equation}\label{eq:rationality}
B=C_{11}+FC_{21}=C_{11}-C_{12}C_{22}^{-1}C_{21}
\end{equation}
(this tells us that $B$ is a rational matrix; thereby giving an alternative proof of Corollary \ref{cor:rat}). 
It is readily checked that    
$$
\begin{pmatrix}
C_{11} & C_{12}\\
C_{21} & C_{22}
\end{pmatrix}
\begin{pmatrix}
I & O \\
-C_{22}^{-1}C_{21} & I
\end{pmatrix}
=\begin{pmatrix}
C_{11}-C_{12}C_{22}^{-1}C_{21}& C_{12}\\
O & C_{22}
\end{pmatrix},
$$
so $\det C=\det (C_{11}-C_{12}C_{22}^{-1}C_{21})\cdot\det C_{22}=\det B\cdot \det C_{22}$. Because $\det C=1$, $\det B={\rm vol}(\mathbb{R}^d/\mathcal{L}_{\mathcal{S}})^2$, and $\det C_{22}={\rm vol}(H_{\mathbb{R}}/H)^2$, our claim follows. \hfill $\Box$

\medskip

The following gives a characterization of crystallographic 1-tight frames by means of a minimal principle.

\begin{prop}\label{prop:minimal}
For a crystallographic frame $\mathcal{S}=\{{\bf v}_i\}_{i=1}^N$ in $\mathbb{R}^d$ with the vanishing group $H$, we have 
$$
\sum_{i=1}^N\|{\bf v}_i\|^2\geq d\cdot{\rm vol}(\mathbb{R}^d/\mathcal{L}_{\mathcal{S}})^{2/d}h(H_{\mathbb{Q}})^{2/d}\big(=d\cdot{\rm vol}(\mathbb{R}^d/\mathcal{L}_{\mathcal{S}})^{2/d}|H^{\#}/H|^{1/d}\big).
$$
The equality holds if and only if $\mathcal{S}$ is tight.
\end{prop}

It should be noted that the quantity $\displaystyle{\rm vol}(\mathbb{R}^d/\mathcal{L}_{\mathcal{S}})^{-2/d}\sum_{i=1}^N\|{\bf v}_i\|^2$ depends only on the similarity class of the crystallographic frame $\mathcal{S}$.

\smallskip

\noindent{\it Proof}.~ For any positive symmetric matrix $B\in M_{d}(\mathbb{R})$, we have 
$$
{\rm tr}~\!B\geq d \big(\det B\big)^{1/d},
$$
where the equality holds if and only if $B=\alpha I_d$ for some $\alpha>0$ (this is easily deduced from the {\it inequality of arithmetic and geometric means}). Applying this inequality to $B={}^t\!AA$, we find that 
$$
\sum_{i=1}^N\|{\bf v}_i\|^2\geq d (\det S_{\mathcal{S}})^{1/d}
$$
for any frame $\mathcal{S}=\{{\bf v}_i\}_{i=1}^N$, 
where the equality holds if and only if $\mathcal{S}$ is tight.

We now let $S$ be the frame operator associated with $\mathcal{S}$. Recall that $S$ is symmetric and positive, whence one can find a (unique) symmetric positive operator $S^{-1/2}$ with $(S^{-1/2})^2=S^{-1}$. Then
\begin{eqnarray*}
&&\sum_{i=1}^N\big\langle {\bf x}, S^{-1/2}{\bf v}_i\big\rangle S^{-1/2}{\bf v}_i
=S^{-1/2}\sum_{i=1}^N\big\langle S^{-1/2}{\bf x},{\bf v}_i\big\rangle {\bf v}_i\\
&=&S^{-1/2}SS^{-1/2}{\bf x}={\bf x}.
\end{eqnarray*}
Therefore $\mathcal{U}=\{S^{-1/2}{\bf u}_i\}_{i=1}^N$ is a 1-tight frame whose vanishing subgroup is obviously $H$. Thus $\rho_{\mathcal{U}}=S^{-1/2}\rho_{\mathcal{S}}$ and 
$$
{\rm vol}\big(\mathbb{R}^d/{\rm Image}~\!\rho_{\mathcal{U}}\big)=\det S^{-1/2}{\rm vol}\big(\mathbb{R}^d/{\rm Image}~\!\rho_{\mathcal{S}}\big),
$$
from which it follows that
$$
\det S={\rm vol}\big(\mathbb{R}^d/\mathcal{L}_{\mathcal{S}}\big)^2{\rm vol}(H_{\mathbb{R}}/H)^2.
$$
This 
completes the proof. 
\hfill $\Box$

\medskip

In Theorem \ref{thm:2dimension}, we observed that a square free positive integer $D$ is associated with each 2-dimensional tight frame $\mathcal{S}$. We close this section with establishing a relationship between the integer $D$ and the height of $H(\mathcal{S})_{\mathbb{Q}}$.

\begin{prop} Let $[z_1,\ldots,z_N]\in Q_N$ 
be a $\mathbb{Q}(\sqrt{-D})$-rational point corresponding to a 2-dimensional tight frame $\mathcal{S}=\{{\bf v}_i\}_{i=1}^N$ whose vanishing group is $H$. 
Then 
$D$ is the square free part of $h(H_{\mathbb{Q}})^2(={\rm vol}(H_{\mathbb{R}}/H)^2$).

\end{prop}

\noindent{\it Proof}.~ One may assume $z_i\in \mathbb{Q}(\sqrt{-D})$. 
Let $w_1,w_2\in \mathbb{C}$ be a $\mathbb{Z}$-basis of the period lattice $\mathcal{L}_{\mathcal{S}}$ (we are working in $\mathbb{C}$ instead of $\mathbb{R}^2$). Writing $w_i=a_i+b_i\sqrt{-D}\in \mathbb{Q}(\sqrt{-D})$, we find 
\begin{eqnarray*}
&&{\rm vol}(\mathbb{C}/\mathcal{L}_{\mathcal{S}})^{-1}\sum_{i=1}^N\|{\bf v}_i\|^2={\rm vol}(\mathbb{C}/\mathcal{L}_{\mathcal{S}})^{-1}(|z_1|^2+\cdots +|z_N|^2)
\\
&& = \frac{1}{|a_1b_2-a_2b_1|\sqrt{D}}(|z_1|^2+\cdots +|z_N|^2),
\end{eqnarray*}
which is equal to $2{\rm vol}(H_{\mathbb{R}}/H)$ in view of proposition \ref{prop:minimal}. Because $|z_1|^2+\cdots +|z_N|^2$ is rational, the claim is proved.   \hfill $\Box$

\section{Tight frames associated with finite graphs and combinatorial Abel's theorem}\label{sect:finite graphs}
The notions of {\it Jacobian} and {\it Picard group} together with {\it Abel's theorem} play a significant role in classical algebraic geometry. 
The aim of this section is to introduce combinatorial analogues of these notions by using certain tight frames associated with finite graphs.

We first fix some basic notations and terminology. A {\it graph}\index{graph} is represented by an ordered pair $X = (V,E)$ of the set of {\it vertices} $V$ and the set of all {\it directed edges} $E$ (note that each edge has just two directions, which are to be expressed by arrows). For an  directed edge $e$, we denote by ${\it o}(e)$ the {\it origin}, and by ${\it t}(e)$ the {\it terminus}. The inversed edge of $e$ is denoted by $\overline{e}$. With these notations, we have $o(\overline{e})=t(e)$, $t(\overline{e})=o(e)$. An {\it orientation} of $X$ is a subset $E^o$ of $E$ such that $E=E^o\cup \overline{E^o}$ and $E^o\cap \overline{E^o}=\emptyset$. We use the notation $E_x$ for the set of directed edges $e$ with $o(e)=x$. Throughout, the degree ${\rm deg}~\!x=|E_x|$ is assumed to be greater than or equal to three for every vertex $x$. 

We let $X_0=(V_0,E_0)$ be a finite connected graph that is regarded as a 1-dimensional cell complex. From now on, we shall make use of (co)homology theory of cell complexes. Let $K$ be $\mathbb{Z}, \mathbb{Q}$, or $\mathbb{R}$, and let $\partial: C_1(X_0,K)\longrightarrow C_0(X_0,K)$ be the boundary operators of chain groups; namely the homomorphism defined by $\partial(e)=t(e)-o(e)$ ($e\in E_0$). The 1-homology group $H_1(X_0,\mathbb{Z})={\rm ker}~\!\partial$ is a direct summand of the 1-chain group $C_1(X_0,\mathbb{Z})$, and is a lattice of $H_1(X_0,\mathbb{R})$. We denote by $b_1=b_1(X_0)$ the {\it betti number}, i.e., $b_1={\rm dim}~\!H_1(X_0,\mathbb{R})$.

Define the natural inner products on $C_0(X_0,\mathbb{R})$ and $C_1(X_0,\mathbb{R})$ by setting
$$
\langle x,y\rangle=\begin{cases}
1 & (x=y)\\
0 & (x\neq y),
\end{cases}
$$
and 
$$
\langle e, e'\rangle=\begin{cases}
1   & (e'=e)\\
-1 &  (e'=\overline{e})\\
0  &  (\text{otherwise}).
\end{cases}
$$
The set of vertices $V_0$ constitutes an orthonormal basis of $C_0(X_0,\mathbb{R})$, while an orientation $E_0^o$ of $X_0$ yields an orthonormal basis of $C_1(X_0,\mathbb{R})$. With the inner product on $H_1(X_0,\mathbb{R})$ induced from the one on $C_1(X_0,\mathbb{R})$, the lattice $H_1(X_0,\mathbb{Z})$ is integral,
and hence the dual lattice $H_1(X_0,\mathbb{Z})^{\#}$
contains $H_1(X_0,\mathbb{Z})$.

Let $\partial^{*}:C_0(X_0,\mathbb{R})\longrightarrow C_1(X_0,\mathbb{R})$ be the adjoint operator of $\partial$ with respect to the above inner products. It is straightforward to see  
$$
\partial^{*} x=-\sum_{e\in E_{0x}}e\quad (x\in V_0),
$$
so $\partial^{*}\big(C_0(X_0,\mathbb{Z})\big)\subset C_1(X_0,\mathbb{Z})$. We also have $H_1(X_0,\mathbb{R})=\big({\rm Image}~\!\partial^{*}\big)^{\perp}$, and hence $C_1(X_0,\mathbb{R})=H_1(X_0,\mathbb{R})\oplus {\rm Image}~\!\partial^{*}$ (orthogonal direct sum).

We denote by $P_0:C_1(X_0,\mathbb{R})\longrightarrow H_1(X_0,\mathbb{R})$ the orthogonal projection, and put ${\bf v}_0(e)=P_0(e)$ ($e\in E_0$). Since ${\rm ker}~\!P_0={\rm Image}~\!\partial^{*}$, we get 
\begin{equation}\label{eq:harmonicityproperty}
\sum_{e\in E_{0x}}{\bf v}_0(e)={\bf 0} \quad (x\in V_0).
\end{equation}

\begin{lemma}\label{lemma:framehomology}
$\mathcal{S}=\big\{{\bf v}_0(e)\big\}_{e\in E_0^o}$ is a crystallographic $1$-tight frame (and hence $\big\{{\bf v}_0(e)\big\}_{e\in E_0}$ is $2$-tight). Its vanishing group is $\partial^{*}\big(C_0(X_0,\mathbb{Z})\big)$, and its period lattice is $H_1(X_0,\mathbb{Z})^{\#}$.
\end{lemma}

\noindent{\it Proof}.~ The first claim is obvious because the orientation $E_0^o$ gives an orthonormal basis of $C_1(X_0,\mathbb{R})$, and ${\bf v}_0(e)$'s are the projected images of this orthonormal basis by the orthogonal projection $P_0$.

For the second claim, it suffices to show that 
$$
({\rm Image}~\!\partial^{*})\cap C_1(X_0,\mathbb{Z})=\partial^{*}\big(C_0(X_0,\mathbb{Z})\big).
$$ 
To this end,
take any 0-chain $\alpha=\displaystyle\sum_{x\in V_0}a_xx\in C_0(X_0,\mathbb{R})$ such that 
$\partial^{*}\alpha\in C_1(X_0,\mathbb{Z})$. 
$$
\partial^{*}\alpha=-\sum_{x\in V_0}\sum_{e\in E_{0x}}a_xe
=-\sum_{e\in E_0}a_{o(e)}e=\sum_{e\in E_0^o}\big(a_{t(e)}-a_{o(e)}\big)e,
$$
so $a_{t(e)}-a_{o(e)}\in \mathbb{Z}$ for every $e\in E_0$. This implies that there exists a real number $a$ with $a_x+a\in \mathbb{Z}$ ($x\in V_0$). Putting 
$$
\beta=\displaystyle\sum_{x\in V_0}(a_x+a)x\in C_0(X_0,\mathbb{Z}),
$$
we obtain $\partial^{*}\alpha=\partial^{*}\beta$. This proves the claim.

To prove that $P_0\big(C_1(X_0,\mathbb{Z})\big)=H_1(X_0,\mathbb{Z})^{\#}$, take a spanning tree $T$ of $X_0$, and let 
$$
e_1,e_2,\ldots,e_{b_1},\overline{e_1},\overline{e_2},\ldots,\overline{e_{b_1}} 
$$ 
be all directed edges not in $T$. The vectors ${\bf v}_0(e_1),\ldots,{\bf v}_0(e_{b_1})$ constitute a $\mathbb{Z}$-basis of the lattice group $H_1(X_0,\mathbb{Z})^{\#}$. This is so because we may create a $\mathbb{Z}$-basis of $H_1(X_0,\mathbb{Z})$ consisting of circuits\index{circuit} $c_1,\ldots,c_{b_1}$ in $X_0$ such that $c_i$ contains $e_i$, and 
$$\langle c_i,{\bf v}_0(e_j)\rangle=\langle c_i,P_0(e_j)\rangle=\langle P_0(c_i),e_j\rangle=\langle c_i,e_j\rangle=\delta_{ij},
$$ 
namely $\big\{{\bf v}_0(e_1),\ldots,{\bf v}_0(e_{b_1})\big\}$ is the dual basis of $\{c_1,\ldots, c_{b_1}\}$. From this, our claim immediately follows.
\hfill $\Box$

\medskip

Let us exhibit two instructive examples.

\smallskip

\begin{figure}[htbp]
\begin{center}
\includegraphics[width=.3\linewidth]
{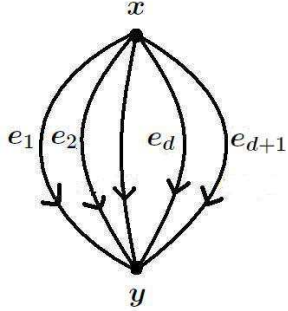}
\end{center}
\caption{Graph $\Delta_{d+1}$}\label{fig:graph}
\end{figure}

(1)~ Let $\Delta_{d+1}$ be the graph depicted in Fig.~\!\ref{fig:graph}. 

Then $\{e_1,\ldots,e_{d+1}\}$ is an orthonormal basis of $C_1(\Delta_{d+1},\mathbb{R})$, and $e_1-e_2$, $e_2-e_3$,\ldots,$e_{d}-e_{d+1}$ is a $\mathbb{Z}$-basis of $H_1(\Delta_{d+1},\mathbb{Z})$. One can check 
$$
{\bf v}_0(e_i)=e_i-\frac{1}{d+1}\sum_{j=1}^{d+1}e_j.
$$ 
Thus $\big\{{\bf v}_0(e_i)\big\}_{i=1}^{d+1}$ is the tight frame associated with the equilateral simplex (see the proof of Proposition \ref{prop:equi}).

It should be pointed out that ${\bf v}_0(e_i)-{\bf v}_0(e_j)=e_i-e_j~(i\neq j)$ are vectors representing ridges (edges) of the simplex, and that $\varPhi=\{e_i-e_j|~ i\neq j\}$ is the irreducible root system of type $A_{d}$ with {\it simple roots} $e_1-e_2$, $e_2-e_3$,\ldots,$e_{d}-e_{d+1}$. 

This example has something to do with the diamond crystal (in the case $d=3$).

\smallskip

\begin{figure}[htbp]
\begin{center}
\includegraphics[width=.3\linewidth]
{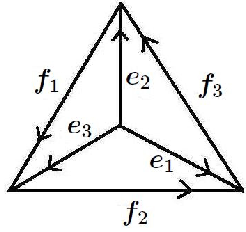}
\end{center}
\caption{Graph $K_4$}\label{fig:k4}
\end{figure}

(2)~ Related to the diamond twin mentioned in Introduction is the complete graph $K_4$ with 4 vertices. See Fig.~\!\ref{fig:k4} for an orientation (the labeling of edges).

Take four closed paths $c_1=(e_2,f_1,\overline{e_3})$, $c_2=(e_3,f_2,\overline{e_1})$, $c_3=(e_1,f_3,\overline{e_2})$, $c_4=(\overline{f_1},\overline{f_2},\overline{f_3})$ in $K_4$. Then $c_1$, $c_2$, $c_3$, $c_4$ as 1-chains 
satisfy $c_1+c_2+c_3+c_4=0$, $\|c_1\|^2=\|c_2\|^2=\|c_3\|^2=\|c_4\|^2$ and $\langle c_i, c_j\rangle=-1$ $~(i\neq j)$ \big($c_1,c_2,c_3$ constitute a $\mathbb{Z}$-basis of $H_1(K_4,\mathbb{Z})$\big). Moreover $c_1$, $c_2$, $c_3$, $c_4$ are vectors represented by directed segments joining the origin and vertices of the regular tetrahedron. An easy computation gives  
\begin{eqnarray*}
&&{\bf v}_0(e_1)=\frac{1}{4}(c_3-c_2),\quad {\bf v}_0(e_2)=\frac{1}{4}(c_1-c_3),\quad  {\bf v}_0(e_3)=\frac{1}{4}(c_2-c_1),\\
&& {\bf v}_0(f_1)=\frac{1}{4}(c_1-c_4),\quad 
 {\bf v}_0(f_2)=\frac{1}{4}(c_2-c_4),\quad
 {\bf v_0}(f_3)=\frac{1}{4}(c_3-c_4).
\end{eqnarray*}
This implies that $\big\{\pm{\bf v}_0(e_i), \pm{\bf v}_0(f_j)\big\}$ is the root system $A_3$.

\medskip

Now let us proceed to a combinatorial analogue of Abel's theorem. The classical Abel's theorem in algebraic geometry gives a link between regular (holomorphic) maps from a complex projective algebraic curve into two kinds of complex tori (called the Jacobian and Picard group, respectively). In the graph-theoretic setting, certain finite abelian groups are to act as complex tori, and the counterparts of the regular map into the Jacobian is defined by using the tight frame given in Lemma \ref{lemma:framehomology}. The approach explained from now is somewaht different from the one in \cite{su4}, \cite{su5} 
(see also \cite{nag}, \cite{Baker} for different ways to introduce the concepts). 

A key role is played by the direct sum $H_1(X_0,\mathbb{Z})+\partial^{*}\big(C_0(X_0,\mathbb{Z})\big)$. Because this is a sublattice of $C_1(X_0,\mathbb{Z})$, the factor group$$
C_1(X_0,\mathbb{Z})/\big[H_1(X_0,\mathbb{Z})+\partial^{*}\big(C_0(X_0,\mathbb{Z})\big)\big]
$$
is a finite abelian group. We shall rewrite this group in two ways; one is $\mathcal{J}(X_0)$, an analogue of Jacobian; another is ${\rm Pic}(X_0)$, an analogue of Picard group.  

In view of the above lemma, the projection $P_0$ induces an isomorphism 
$$
C_1(X_0,\mathbb{Z})/\partial^{*}\big(C_0(X_0,\mathbb{Z})\big)\longrightarrow H_1(X_0,\mathbb{Z})^{\#}.
$$
Notice that the inclusion 
$$
H_1(X_0,\mathbb{Z})\longrightarrow C_1(X_0,\mathbb{Z})
$$ 
induces an injection $H_1(X_0,\mathbb{Z})\longrightarrow C_1(X_0,\mathbb{Z})/\partial^{*}\big(C_0(X_0,\mathbb{Z})\big)$, so that $H_1(X_0,\mathbb{Z})$ is identified with a subgroup of $C_1(X_0,\mathbb{Z})/\partial^{*}\big(C_0(X_0,\mathbb{Z})\big)$. The double factor group 
$$
\big[C_1(X_0,\mathbb{Z})/\partial^{*}\big(C_0(X_0,\mathbb{Z})\big)\big]/H_1(X_0,\mathbb{Z})
$$ 
is isomorphic to $C_1(X_0,\mathbb{Z})/\big[H_1(X_0,\mathbb{Z})+\partial^{*}\big(C_0(X_0,\mathbb{Z})\big)\big]$. Thus we have an isomorphism
$$
C_1(X_0,\mathbb{Z})/\big[H_1(X_0,\mathbb{Z})+\partial^{*}\big(C_0(X_0,\mathbb{Z})\big)\big]\longrightarrow H_1(X_0,\mathbb{Z})^{\#}/H_1(X_0,\mathbb{Z}),
$$
where we should note that $\mathcal{J}(X_0)=H_1(X_0,\mathbb{Z})^{\#}/H_1(X_0,\mathbb{Z})$ is a finite subgroup of the torus group $J(X_0)=H_1(X_0,\mathbb{R})/H_1(X_0,\mathbb{Z})$. One may think of $\mathcal{J}(X_0)$ as a combinatorial analogue of Jacobian. 

On the other hand, the boundary operator induces an isomorphism
$$
C_1(X_0,\mathbb{Z})/H_1(X_0,\mathbb{Z}) \longrightarrow {\rm ker}~\!\epsilon,
$$
where $\epsilon :C_0(X_0,\mathbb{Z})\longrightarrow \mathbb{Z}$ is the homomorphism (called {\it argumentation}) defined by 
$$
\epsilon\Big(\displaystyle\sum_{x\in V_0}a_xx\Big)=\sum_{x\in V_0}a_x.
$$ 
Thus in a similar way as above, we obtain an isomorphism
$$
C_1(X_0,\mathbb{Z})/\big[H_1(X_0,\mathbb{Z})+\partial^{*}\big(C_0(X_0,\mathbb{Z})\big)\big]\longrightarrow {\rm ker}~\!\epsilon/\partial \partial^{*}\big(C_0(X_0,\mathbb{Z})\big),
$$
where ${\rm ker}~\!\epsilon$ is regarded as an analogue of the group of divisors with degree $0$, so that we denote it by ${\rm Div}^0(X_0)$ after the notation in algebraic geometry. Meanwhile $\partial \partial^{*}\big(C_0(X_0,\mathbb{Z})\big)$ is regarded as an analogue of the group of principal divisors (see \cite{Baker} for the reason). Therefore we denote it by ${\rm Prin}(X_0)$. Remembering again the terminology in algebraic geometry, it is justified to write ${\rm Pic}(X_0)$ for the factor group ${\rm Div}^0(X_0)/{\rm Prin}(X_0)$ and to call it the Picard group of $X_0$.

Summarizing the argument, we have

\begin{prop}
There is a natural isomorphism $\varphi$ of the Picard group ${\rm Pic}(X_0)$ onto the Jacobian $\mathcal{J}(X_0)$. 
\end{prop}

To describe $\varphi$ more explicitly, pick up a reference vertex $x_0\in V_0$.  For $x\in V_0$, select a path $c=(e_1,\ldots,e_n)$ in $X_0$ such that $o(c)(=o\big(e_1)\big)=x_0$ and $t(c)\big(=t(e_n)\big)=x$. Then regarding $c$ as a 1-chain, we get 
$\partial c=x-x_0$. 
From the way to construct $\varphi$, we easily find that $\varphi$ brings $x-x_0\in {\rm ker}~\!\epsilon$ modulo $\partial \partial^{*}\big(C_0(X_0,\mathbb{Z})\big)$ to ${\bf v}_0(e_1)+\cdots+{\bf v}_0(e_n)\in H_1(X_0,\mathbb{Z})^{\#}$ modulo $H_1(X_0,\mathbb{Z})$. 

Now define the {\it combinatorial Albanese map} $\varPhi^{\rm al}:V\longrightarrow \mathcal{J}(X_0)$ by setting
$$
\varPhi^{\rm al}(x)={\bf v}_0(e_1)+\cdots+{\bf v}_0(e_n)\in H_1(X_0,\mathbb{Z})^{\#}~\text{modulo}~ H_1(X_0,\mathbb{Z}),
$$
where the sum on the right hand side is an analogue of the line integral of a holomorphic 1-form which appears in the definition of classical Albanese maps (we shall see in the next section that ${\bf v}_0$ deserves to be called ``harmonic" as a cochain of $X_0$). 
We also define the {\it combinatorial Abel-Jacobi map} $\varPhi^{\rm aj}:V_0\longrightarrow {\rm Pic}(X_0)$ by
$$
\varPhi^{\rm aj}(x)=x-x_0\in {\rm Div}^0(X_0)~ \text{modulo}~{\rm Prin}(X_0).
$$
By definition, we obtain $\varPhi^{\rm al}=\varphi \circ \varPhi^{\rm aj}$. This is nothing but an analogue of classical Abel's theorem.

To describe the structure of $\mathcal{J}(X_0)$, we consider the integral matrix $A=\big(\langle\alpha_i,\alpha_j\rangle\big)$ $\in M_{b_1}(\mathbb{Z})$ where $\{\alpha_1,\ldots,\alpha_{b_1}\}$ is a $\mathbb{Z}$-basis of $H_1(X_0,\mathbb{Z})$. Applying the theory of elementary divisors to $A$, we find $P,Q\in GL_{b_1}(\mathbb{Z})$ and positive integers $k_1,\ldots, k_{b_1}$ such that 
\begin{equation}\label{eq:integralmatrix}
PAQ=\begin{pmatrix}
k_1 & 0 & 0 & \cdots & 0 \\
0   & k_2 & 0 &\cdots & 0 \\
 & \cdots & \cdots & \\
0 & 0 & 0& \cdots & k_{b_1}
\end{pmatrix},
\end{equation}
where $k_i$ divides $k_{i+1}$ ($i=1,\ldots, b_1-1$). The array $(k_1,\ldots,k_{b_1})$, which depends only on $X_0$, determines the structure of $\mathcal{J}(X_0)$; that is, $\mathcal{J}(X_0)=\mathbb{Z}_{k_1}\times\cdots\times \mathbb{Z}_{k_{b_1}}$. For instance, using this fact, we find $\mathcal{J}(\Delta_{d+1})=\mathbb{Z}_{d+1}$, and $\mathcal{J}(K_n)=(\mathbb{Z}_n)^{n-2}$. 

We have more about $\mathcal{J}(X_0)$.
Algebraic graph theory \cite{biggs1}, \cite{biggs2} allows us to establish the fact that the order of $\mathcal{J}(X_0)$ \big(and ${\rm Pic}(X_0)$\big) is equal to $\kappa(X_0)$, the {\it tree number} for $X_0$, which is defined to be the number of spanning trees in $X_0$; therefore $k_1\cdots k_{b_1}=\kappa(X_0)$ (see \cite{su5}). Further the canonical inner product on $H_1(X_0,\mathbb{R})$ induces a flat metric on the torus $J(X_0)=H_1(X_0,\mathbb{R})/H_1(X_0,\mathbb{Z})$ for which, in view of (\ref{eq:dualint}), we have 
$$
{\rm vol}\big(J(X_0)\big)
=\kappa(X_0)^{1/2}.
$$

 The Jacobian $\mathcal{J}(X_0)$ has another appendage. A non-degenerate symmetric bilinear form 
on $\mathcal{J}(X_0)$ with values in $\mathbb{Q}/\mathbb{Z}$ is induced from the inner product on $H_1(X_0,\mathbb{R})$. 
One may think of this bilinear form as an analogue of ``principal polarization"

Before closing this section, we make a minor remark on direct summands of $\mathbb{Z}^N$ which is related to the above discussion. 
Let $H$ be a direct summand of $\mathbb{Z}^N$ of rank $k$.
We put
$$
H^{{\perp}_{\mathbb{Z}}}=\{{\bf x}\in \mathbb{Z}^N|~\langle {\bf x},H\rangle=0\},
$$
which is clearly a direct summand of $\mathbb{Z}^N$ of rank $N-k$. We easily observe that $\big(H^{{\perp}_{\mathbb{Z}}}\big)^{{\perp}_{\mathbb{Z}}}=H$.
Although $\big(H^{{\perp}_{\mathbb{Z}}}\big)_{\mathbb{R}}=H_{\mathbb{R}}^{\perp}$ and $H_{\mathbb{R}}\oplus H_{\mathbb{R}}^{\perp}=\mathbb{R}^N$, the direct sum $H\oplus H^{{\perp}_{\mathbb{Z}}}$ does not agree with $\mathbb{Z}^N$ in general. Indeed 
$
\mathbb{Z}^N/\big(H\oplus H^{{\perp}_{\mathbb{Z}}}\big)
$
is isomorphic to $H^{\#}/H$ (this is actually a generalization of what we have seen above). To prove this, let $P:\mathbb{R}^N\longrightarrow \mathbb{R}^N$ be the orthogonal projection whose image is $H_{\mathbb{R}}$. It suffices to check that $P(\mathbb{Z}^N)=H^{\#}$; for if this is true, then $P|\mathbb{Z}^N$ induces an isomorphism of $\mathbb{Z}^N/H^{{\perp}_{\mathbb{Z}}}$ onto $H^{\#}$ because ${\rm Ker}~\!P|\mathbb{Z}^N=H^{{\perp}_{\mathbb{Z}}}$. 
Take a $\mathbb{Z}$-basis $\{{\bf a}_1,\ldots,{\bf a}_N\}$ of $\mathbb{Z}^N$ such that $\{{\bf a}_1,\ldots,{\bf a}_k\}$ is a $\mathbb{Z}$-basis of $H$. Since $A=({\bf a}_1,\ldots,{\bf a}_N)\in GL_N(\mathbb{Z})$, there exists $B=({\bf b}_1,\ldots,{\bf b}_N)\in GL_N(\mathbb{Z})$ such that ${}^tBA=I_N$, or equivalently $\langle {\bf b}_i,{\bf a}_j\rangle=\delta_{ij}$. Therefore 
$$
\big\langle P({\bf b}_i),{\bf a}_j\big\rangle=
\big\langle {\bf b}_i,P({\bf a}_j)\big\rangle=
\langle {\bf b}_i,{\bf a}_j\rangle=\delta_{ij}\quad (1\leq j\leq k),
$$
which implies  
$\big\{P({\bf b}_1),\ldots,P({\bf b}_k)\big\}$ is a $\mathbb{Z}$-basis of $H^{\#}$.
Therefore $P(\mathbb{Z}^N)=H^{\#}$, as required.

Using $\big|\mathbb{Z}^N/\big(H\oplus H^{{\perp}_{\mathbb{Z}}}\big)\big|=|H^{\#}/H|$, we find
\begin{eqnarray*}
{\rm vol}(H_{\mathbb{R}}/H){\rm vol}(H_{\mathbb{R}}^{\perp}/H^{{\perp}_{\mathbb{Z}}})&=&{\rm vol}\big(\mathbb{R}^N/(H\oplus H^{{\perp}_{\mathbb{Z}}})\big)\\
&=&\big|\mathbb{Z}^N/\big(H\oplus H^{{\perp}_{\mathbb{Z}}}\big)\big|{\rm vol}(\mathbb{R}^N/\mathbb{Z}^N)\\
&=& |H^{\#}/H|,
\end{eqnarray*}
so 
$$
{\rm vol}(H_{\mathbb{R}}^{\perp}/H^{{\perp}_{\mathbb{Z}}})={\rm vol}(H_{\mathbb{R}}/H)=|H^{\#}/H|^{1/2}.
$$

\section{Standard crystal models and tight frames}\label{sect:Standard crystal models}

The combinatorial Albanese map $\varPhi^{\rm al}:V_0\longrightarrow \mathcal{J}(X_0)$ extends to a piecewise linear map of $X_0$ into the flat torus $J(X_0)=H_1(X_0,\mathbb{R})/H_1(X_0,\mathbb{Z})$:
$$
\varPhi^{\rm al}:X_0\longrightarrow J(X_0),
$$
which, if we think of $X_0$ as a (singular) Riemannian manifold, turns out to be {\it harmonic} in the sense of Eells and Sampson \cite{es} (see \cite{sk2} for the detail).

Let $X_0^{\rm ab}$ be the maximal abelian covering graph over $X_0$, i.e., the abelian covering graph over $X_0$ whose covering transformation group is $H_1(X_0,\mathbb{Z})$. Consider a lifting $\varPhi_0:X_0^{\rm ab}\longrightarrow H_1(X_0,\mathbb{R})$ of $\varPhi^{\rm al}$, which obviously satisfies
$$
\varPhi_0(\sigma x)=\varPhi_0(x)+\sigma \quad \big(\sigma\in H_1(X_0,\mathbb{Z})\big).
$$
The image $\varPhi_0(X_0^{\rm ab})$ is considered a $b(X_0)$-dimensional crystal net. For instance, if $X_0=\Delta_4$ (resp. $X_0=K_4$), then $\varPhi_0(X_0^{\rm ab})$ is the diamond crystal (resp. the diamond twin); see the previous section. 
 
Having this observation in mind, we consider general abelian covering graphs over $X_0$ and their realizations. A $d$-dimensional {\it topological graph} is an infinite-fold abelian covering graph $X=(V,E)$ over a finite graph $X_0$ whose covering transformation group is a free abelian group $L$ of rank $d$. Theory of covering spaces tells us that there is a subgroup $H$ (called a {\it vanishing subgroup}) such that $H_1(X_0,\mathbb{Z})/H=L$ \big(note that $H$ is a direct summand of $H_1(X_0,\mathbb{Z})$\big). Actually the topological crystal $X$ is the quotient graph of $X_0^{\rm ab}$ over $X_0$ modulo $H$. In this view, we call $X_0^{\rm ab}$ the {\it maximal topological crystal} over $X_0$.

A (periodic) {\it realization} is a piecewise linear map $\varPhi:X\longrightarrow \mathbb{R}^d$ satisfying 
$$
\varPhi(\sigma x)=\varPhi(x)+\rho(\sigma)\qquad (\sigma \in L),
$$
where $\rho:L\longrightarrow \mathbb{R}^d$ is an injective homomorphism whose image is a lattice in $\mathbb{R}^d$.\footnote{The network $\varPhi(X)$ could be ``degenerate" in the sense that different vertices of $X$ are realized as one points, or different edges overlap in $\mathbb{R}^d$. But we shall not exclude these possibilities.} We call $\rho$ \big(resp. $\rho(L)$\big) the {\it period homomorphism} (resp. the {\it period lattice}) for $\varPhi$.

By putting ${\bf v}(e)=\varPhi\big(t(e)\big)-\varPhi\big(o(e)\big)$ ~$(e\in E)$, we obtain a $L$-invariant function ${\bf v}$ on $E$ which we may identify with a 1-cochain ${\bf v}\in C^1(X_0,\mathbb{R}^d)$ with values in $\mathbb{R}^d$. Since ${\bf v}$ determines completely $\varPhi$ (up to parallel translations), we shall call ${\bf v}$ the {\it building cochain}\footnote{In \cite{su4}, \cite{su5}, the term ``building block" is used. The idea to describe crystal structures by using finite graphs together with vector labeling is due to \cite{chung}} of $\varPhi$. One can check that if we identify the cohomology class $[{\bf v}]\in H^1(X_0,\mathbb{R}^d)$ with a homomorphism of $H_1(X_0,\mathbb{Z})$ into $\mathbb{R}^d$ (the {\it duality of cohomology and homology}), then $[{\bf v}]=\rho\circ \mu$, where $\mu:H_1(X_0,\mathbb{Z})\longrightarrow L$ is the canonical homomorphism. In particular, ${\rm Ker}~\![{\bf v}]=H$ and ${\rm Image}~\![{\bf v}]=\rho(L)$. 

\begin{lemma}{\rm \cite{su4}}
Giving a periodic realization of a topological crystal over $X_0$ is equivalent to giving a 1-cochain ${\bf v}\in C^1(X_0,\mathbb{R}^d)$ such that the image of the homomorphism $[{\bf v}]:H_1(X_0,\mathbb{Z})\longrightarrow \mathbb{R}^d$ is a lattice in $\mathbb{R}^d$. 
\end{lemma}

We are now at the stage to give the definition of standard realizations. 
Let $H_{\mathbb{R}}$ be the subspace of $H_1(X_0,\mathbb{R})$ spanned by the vanishing group $H$, 
and $H_{\mathbb{R}}^{\perp}$ the orthogonal complement of $H_{\mathbb{R}}$ in $H_1(X_0,\mathbb{R})$:
$$
H_1(X_0,\mathbb{R})=H_{\mathbb{R}}\oplus H_{\mathbb{R}}^{\perp}.
$$
Then ${\rm dim}~\hspace{-0.05cm}H_{\mathbb{R}}^{\perp}={\rm rank}~\hspace{-0.05cm}L=d$. By choosing an orthonormal basis of $H_{\mathbb{R}}^{\perp}$, we identify $H_{\mathbb{R}}^{\perp}$ with the Euclidean space $\mathbb{R}^d$. Using the orthogonal projection $P:H_1(X_0,\mathbb{R})\longrightarrow H_{\mathbb{R}}^{\perp}$, we put ${\bf v}(e)=P\big({\bf v}_0(e)\big)$. Then one can check that ${\bf v}$ is the building cochain of a realization $\varPhi:X\longrightarrow \mathbb{R}^d$. We call $\varPhi$ the {\it normalized standard realization} of $X$. If we say simply ``standard realization", it means a realization obtained by performing a similar transformation to the normalized one.

In view of the properties of ${\bf v}_0$ established in the previous section, we find

\smallskip

(1) ({\bf Harmonicity})
\begin{equation}\label{eq:harmonic}
\displaystyle\sum_{e\in E_{0x}}{\bf v}(e)={\bf 0} \quad (x\in V_0),
\end{equation}

\smallskip
(2) ({\bf Tight-frame condition}) 
\begin{equation}\label{eq:standard}
\sum_{e\in E_0} \langle{\bf x},{\bf v}(e) \rangle {\bf v}(e)=2{\bf x}\quad ({\bf x}\in \mathbb{R}^d).
\end{equation}

Furthermore $\big\{{\bf v}(e)\big\}_{e\in E_0}$ is crystallographic; namely it generates a lattice.  Since this lattice contains the image of $[{\bf v}]$, the period lattice is essentially rational. 

\begin{remark}
The direct sum $H\oplus \partial^{*}\big(C_0(X_0,\mathbb{Z})\big) \big(\subset C_1(X_0,\mathbb{Z})\big)$ is contained in the vanishing subgroup for the tight frame $\{{\bf v}(e)\}_{e\in E_0^o}$; however they need not coincide with each other.

\end{remark}

\begin{figure}[htbp]
\begin{center}
\includegraphics[width=.9\linewidth]
{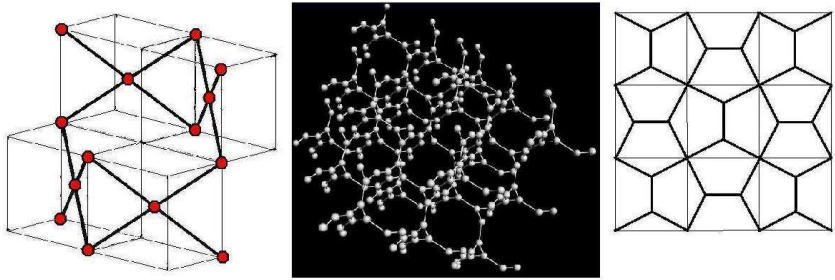}
\end{center}
\caption{Standard realizatins}\label{fig:more}
\end{figure}

Figure \ref{fig:more} exhibits a few more examples of standard realizations (the picure on the right side is a tiling of pentagons with picturesque properties that has become known as the {\it Cairo pentagon}).

Remarkable features of standard realizations are embodied in 

\begin{prop}\label{prop:maximal} {\rm (1)} The standard realization of $X$ is the unique minimizer, up to similar transformations, of the energy\footnote{The energy defined here is considered the potential energy per unit cell when we think of the crystal net as a system of {\it harmonic oscillators}. Clearly this is similarity-invariant.} 
$$
\mathcal{E}(\varPhi)={\rm vol}\big(\mathbb{R}^d/\rho(L)\big)^{-2/d}\sum_{e\in E_0}\|{\bf v}(e)\|^2.
$$

{\rm (2)} Let $\varPhi:X\longrightarrow \mathbb{R}^d$ be the standard realization. Then there exists a homomorphism $\kappa$ of the automorphism group ${\rm Aut}(X)$ of $X$ into the congruence group $M(d)$ of $\mathbb{R}^d$ such that

{\rm (a)} when we write $\kappa(g)=\big(A(g),b(g)\big)\in O(d)\times \mathbb{R}^d$, we have
$$
\varPhi(gx)=A(g)\varPhi(x)+b(g) \qquad (x\in V), 
$$

{\rm (b)} the image $\kappa\big({\rm Aut(X)}\big)$ is a crystallographic group\index{crystallographic group}, a discrete co-compact subgroup of the motion group of $\mathbb{R}^d$ $(see $\cite{charlap}$)$.
\end{prop}

The assertion (1) is a direct consequence of Proposition \ref{prop:minimal}. The second one tells us that the standard realization has maximal symmetry. See \cite{su5} for the proof relying on an asymptotic property of random walks on topological crystals (see also \cite{sk1},\cite{sk2},\cite{sk11}).

\medskip

Equation (\ref{eq:harmonic}) says that the cochain ${\bf v}$ is ``harmonic" in the sense that $\delta{\bf v}=0$ where $\delta:C^1(X_0,\mathbb{R}^d)\longrightarrow C^1(X_0,\mathbb{R}^d)$ is the adjoint of the coboundary operator $d:C^0(X_0,\mathbb{R}^d)\longrightarrow C^1(X_0,\mathbb{R}^d)$ with respect to the natural inner products in $C^i(X_0,\mathbb{R}^d)$. A discrete analogue of the Hodge--Kodaira theorem, which is almost trivial, assures us that the correspondence ${\bf v}\in {\rm Ker}~\!\delta\mapsto [{\bf v}]\in H^1(X_0,\mathbb{R}^d)$ is a linear isomorphism \big(hence ${\rm dim}~\!{\rm Ker}~\!\delta=db_1(X_0)$\big). Thus given $\rho$, there is a unique harmonic cochain ${\bf v}$ with $[{\bf v}]=\rho\circ \mu$. A realization satisfying (\ref{eq:harmonic}) is said to be a {\it harmonic realization} \cite{sk2} (or an equilibrium placement \cite{d-1}), which is characterized as a minimizer of $\mathcal{E}$ when $\rho$ is fixed\footnote{Looking at things through {\it discrete geometric analysis} \cite{su3}, one can see lots of conceptual resemblance between crystallography and electric circuits.}.

Let ${\bf v}\in C^1(X_0,\mathbb{R}^d)$ be a general building cochain satisfying  $[{\bf v}]=\rho\circ \mu$.
The {\it distortion} of the realization given by ${\bf v}$ from the harmonic one is measured by the 0-cochain ${\bf f}\in C^0(X_0,\mathbb{R}^d)$ defined by
\begin{equation}\label{eq:resultant}
{\bf f}(x)=\sum_{e\in E_{0x}}{\bf v}(e)\big(=-(\delta{\bf v})(x)\big),
\end{equation}
which is considered the resultant force acting on the ``atom" $x$ when we regard the crystal net as a system of harmonic oscillators. Obviously 
\begin{equation}\label{eq:force}
\displaystyle\sum_{x\in V_0}{\bf f}(x)={\bf 0}.
\end{equation}
Conversely if ${\bf f}\in C^0(X_0,\mathbb{R}^d)$ satisfies (\ref{eq:force}), then there exists a unique building cochain ${\bf v}$ satisfying (\ref{eq:resultant}) and $[{\bf v}]=\rho\circ \mu$. Indeed this is a consequence of ${\rm Image}~\!\delta=({\rm ker}~\!d)^{\perp}$ and ${\rm ker}~\!\delta=({\rm Image}~\!d)^{\perp}$.

Eigenvalues of the frame operator $S:\mathbb{R}^d\longrightarrow\mathbb{R}^d$ associated with a building cochain ${\bf v}$ gives an information about how much the harmonic realization $\varPhi$ is distorted from the standard one. More precisely, if $\lambda^{\min}$ (resp. $\lambda^{\max}$) is the minimal (resp. maximal) eigenvalue of $S$, then the ratio $R(\varPhi)=\lambda^{\max}/\lambda^{\min}(\geq 1)$ 
is considered representing the degree of distorsion. Indeed $R(\varPhi)=1$ if and only if $\varPhi$ is standard. 

We associate the flat torus $H_{\mathbb{R}}/H$ with a vanishing subgroup $H\subset H_1(X_0,\mathbb{Z})$. As before, we easily observe that ${\rm vol}(H_{\mathbb{R}}/H)^2$ is an integer, and can prove, by modifying slightly the argument used in the proof of Proposition \ref{prop:heightfunction}, that for any $c>0$, there are only finitely many $H$ of rank $b_1-d$ such that ${\rm vol}(H_{\mathbb{R}}/H)<c$. For this, we just work in $H_1(X_0,\mathbb{Z})$ instead of $\mathbb{Z}^N$.

Let $\mathcal{L}_H~(\subset \mathbb{R}^d)$ be the period lattice for the normalized standard realization of the topological crystal corresponding to $H$, and put $J(X_0,H)=\mathbb{R}^d/\mathcal{L}_H$ \big(thus $J(X_0,\{0\})=J(X_0)$\big). Imitating the proof of Proposition \ref{prop:heightvolume}, one may prove 
$$
{\rm vol}\big(J(X_0,H)\big)={\rm vol}\big(J(X_0)\big)/{\rm vol}(H_{\mathbb{R}}/H)=\kappa(X_0)^{1/2}/{\rm vol}(H_{\mathbb{R}}/H).
$$

As an application of this fact, we take up the issue of ``reality" of the standard realization; namely we ask how much part of the family of crystal models is occupied by standard ones which look like genuine crystals. A rough answer is that if we fix the base graph $X_0$, ``most" standard realizations do not look realistic. 

To be more precise, we start with the inequality
$$
\sum_{e\in E_0}\|{\bf v}(e)\|^2\|{\bf x}\|^2\geq \sum_{e\in E_0}\langle{\bf v}(e),{\bf x} \rangle^2=2\|{\bf x}\|^2,
$$
from which we get 
\begin{equation}\label{eq:realistic2}
\max_{e\in E_0}\|{\bf v}(e)\|\geq (2/|E_0|)^{1/2}.
\end{equation}
This implies that the maximal length of ${\bf v}(e)$ is bounded from below by a positive constant depending only on the base graph $X_0$. 
On the other hand, there exists a positive constant $c_d$ such that for any lattice group $\mathcal{L}$ in $\mathbb{R}^d$
$$
\min_{{\bf x}\in \mathcal{L}\backslash\{{\bf 0}\}}\|{\bf x}\|\leq c_d{\rm vol}(\mathbb{R}^d/\mathcal{L})^{1/d}
$$
(this is the celebrated ``convex body theorem" due to Minkowski; cf. \cite{milnor}). 
Applying this fact to the lattice group $\mathcal{L}_H$, we find a non-zero ${\bf x}=\rho(\sigma)\in \mathcal{L}_H$ such that $\|\rho(\sigma)\|\leq c_d{\rm vol}\big(J(X_0,H)\big)^{1/d}$, and hence 
\begin{equation}\label{eq:realistic1}
\|\varPhi(\sigma x)-\varPhi(x)\|\leq c_d{\rm vol}\big(J(X_0,H)\big)^{1/d}.
\end{equation} 

These facts enable us to establish the following theorem.

\begin{thm}
Let $c$ be a positive constant with $c\leq 1$. For a fixed $X_0$, there are only finitely many 
$d$-dimensional topological crystals $X$ over $X_0$ whose standard realization satisfy 

\smallskip

{\rm (i)} $\varPhi:X\longrightarrow \mathbb{R}^d$ is injective,

\smallskip

{\rm (ii)} $\displaystyle\|\varPhi(x)-\varPhi(y)\|\leq c\max_{e\in E_x}\|{\bf v}(e)\|~\Longrightarrow~ \text{$y$ is adjacent to $x$, or $y=x$}$.

\end{thm}

From the nature of genuine crystals, the first condition sounds natural.  
The second condition roughly means that two ``atoms" close enough to each other must be joined by a bond\footnote{This condition may not be enough (or may be too strong) to characterize the ``reality" of a crytstal model because the physical and chemical aspects of crystals are ignored. In particular, the {\it electron clouds}, which are responsible for chemical bonding in crystals, are not involved in the simple network models.}. It is likely that the conclusion in the theorem is true for any reasonable
definition of ``reality" of crystals.

The proof goes as follows. Suppose that there exist infinitely many $X_n$ satisfying (i) and (ii), and let $H_n$ be the vanishing subgroup corresponding to $X_n$. Then ${\rm vol}\big((H_n)_{\mathbb{R}}/H_n\big)$ tends to infinity, so that ${\rm vol}\big(J(X_0,H_n)\big)$ tends to zero as $n\to \infty$. Take an integer $k$ with $k>\displaystyle\max_{x\in X_0} {\rm deg}~\!x=\max_{x\in X} {\rm deg}~\!x$. For a given $\epsilon>0$ with $\epsilon <c(2/|E_0|)^{1/2}$, choose $n$ such that $c_d{\rm vol}\big(J(X_0,H_n)\big)^{1/d}<\epsilon/k$. By (\ref{eq:realistic1}), one can find a non-zero $\sigma\in L_n=H_1(X_0,\mathbb{Z})/H_n$ such that  $\|\varPhi_n(\sigma^ix)-\varPhi_n(x)\|<\epsilon$, where $\varPhi_n$ is the normalized standard realization of $X_n$. Picking up a vertex $x$ satisfying $\max_{e\in E_x}\|{\bf v}(e)\|\geq (2/|E_0|)^{1/2}$ (see (\ref{eq:realistic2})), and putting $x_i=\sigma^ix$, we obtain $k$ distinct vertices $x_1,\ldots, x_k$ such that$$
\|\varPhi(x_i)-\varPhi(x)\|<\epsilon<c\displaystyle\max_{e\in E_x}\|{\bf v}(e)\|.$$ 
Therefore by the condition (ii), $x_i$'s are adjacent to $x$. This implies thet ${\rm deg}~\!x\geq k$, thereby a contradiction.

\medskip
The set of similarity classes of standard realizations of all $d$-dimensional topological crystals over $X_0$ is identified with the rational Grassmannian ${\rm Gr}_{b_1-d}\big(H_1(X_0,\mathbb{Q})\big)$ \big($b_1=b_1(X_0)$\big). As expected from the discussion in Sect.~\!\ref{sec:Parameterizations of crystallographic tight frames}, there is a special feature of the parameterization of 2D standard realizations. To explain this,
we introduce the complex vector space 
$$
\mathbb{H}=\Big\{{\bf z}\in C^1(X_0,\mathbb{C})|~\sum_{e\in E_{0x}}{\bf z}(e)=0~ (x\in V_0)\Big\}.
$$
This is nothing but the space of harmonic cochains (we are identifying $\mathbb{R}^2$ with $\mathbb{C}$), so we find that $
{\rm dim}_{\mathbb{C}}\mathbb{H}=b_1(X_0)
$
We denote by $P(\mathbb{H})$ the projective space associated with $\mathbb{H}$, and by $Q(X_0)$ 
the quadric defined by
$$
Q(X_0)=\Big\{
[{\bf z}]\in P(\mathbb{H})\big|~\sum_{e\in E_0}{\bf z}(e)^2=0
\Big\}.
$$
If we fix an orientation $E_0^o=\{e_1,\ldots,e_N\}$, then by the correspondence ${\bf z}\mapsto \big({\bf z}(e_1),\ldots,{\bf z}(e_N)\big)$, we may think of $Q(X_0)$ as a $(b_1(X_0)-2)$-dimensional subvariety of the complex projective space $P^{N-1}(\mathbb{C})$. Then the intersection 
$$
Q(X_0)\cap \bigcup_{D>0\atop square free}P^{N-1}\big(\mathbb{Q}(\sqrt{-D})\big)
$$ 
is identified with the family of all oriented similarity classes of standard realizations of 2-dimensional topological crystals over $X_0$ \cite{su6}. This is a straightforward generalization of the observation in Sect.~\!\ref{sec:Parameterizations of crystallographic tight frames}. We may also prove that $D$ is the square free part of $\kappa(X_0){\rm vol}(H_{\mathbb{R}}/H)^2$ provided that $H$ is the vanishing subgroup for the standard realization corresponding to a point in $Q(X_0)\cap P^{N-1}\big(\mathbb{Q}(\sqrt{-D})\big)$.  

\begin{example}
(1)~ The {\it kagome lattice} (Fig.~\!\ref{fig:kagomelattice}) corresponds to the $\mathbb{Q}(\sqrt{-3})$-rational points 
$$
\Big[\frac{1\pm\sqrt{-3})}{2}, \frac{1\mp\sqrt{-3}}{2}, -1, \frac{1\pm\sqrt{-3}}{2}, \frac{1\mp\sqrt{-3}}{2}, -1\Big]
$$
of the 2-dimensional projective quadric
\begin{eqnarray*}
&&\{[z_1,z_2,z_3, z_4,z_5,z_6]\in P^5(\mathbb{C});~z_1{}^2+\cdots+z_6{}^2=0,\\
&& \qquad z_1+z_6=z_3+z_4,~ z_2+z_4=z_1+z_5,~z_3+z_5=z_2+z_6\}.
\end{eqnarray*}

\begin{figure}[htbp]
\begin{center}
\includegraphics[width=.5\linewidth]
{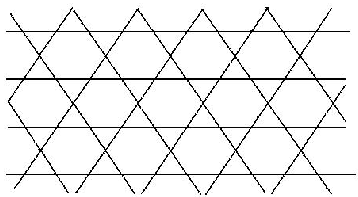}
\end{center}
\caption{Kagome lattice}\label{fig:kagomelattice}
\end{figure}
\end{example}

\begin{figure}[htbp]
\begin{center}
\includegraphics[width=.4\linewidth]
{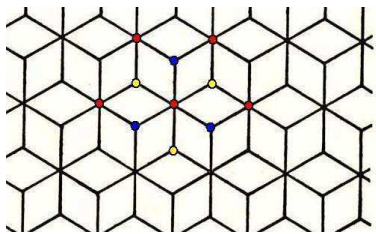}
\end{center}
\caption{Dice lattice}\label{fig:dice}
\end{figure}

(2)~ Figure \ref{fig:dice} is the so-called {\it dice lattice} (also referred to as the $\mathcal{T}_3$ lattice). This corresponds to $\mathbb{Q}(\sqrt{-3})$-rational points
$$
\Big[1, \frac{-1\pm\sqrt{-3}}{2}, \frac{-1\mp\sqrt{-3}}{2}, -1, \frac{1\pm\sqrt{-3}}{2}, \frac{1\mp\sqrt{-3}}{2}\Big]
$$
of the quadric 
$$\{[z_1,\ldots,z_6]\in P^5(\mathbb{C})|~z_1{}^2+\cdots+z_6{}^2=0,~
z_1+z_2+z_3=0,~z_4+z_5+z_6=0\}.
$$

\begin{remark} Two examples above are realizations of 2D topological crystals which come from periodic tilings (tessellations 
of tiles in the plane which are periodic with respect to the translational
action by lattice groups). Easy topological considerations lead to the fact that there are only finitely many 2D topological crystals over a fixed finite graph 
$X_0$ whose standard realizations are 1-skeletons of tilings. In other words, there are only
finitely many rational points in $Q(X_0)$ which correspond to tilings. What we need to notice is that a tiling induces a cellular decomposition of the 2D torus, and that there are only finitely many ways (in a topological sense) to attach 2-cells to a finite graph in order to obtain a torus.
\end{remark}

Proposition \ref{prop:maximal} somehow claims that the standard realization is a natural concept\footnote{The special features of standard realizations might remind the reader of the claim about the {\it golden ratio} ($1.618033\cdots$), a root of $x^2=x+1$, which overemphasizes its significance in the history of art, architecture, sculpture and anatomy. I am not going to overrate the significance of standard realizations though quite a few structures in nature and art are explained using standard realizations.}. To give another justification for the adjective ``standard", we shall see that a crystal net $\varPhi(X)$ with ``big" symmetry is a standard model. 

We assume for simplicity that 
$\varPhi:X\longrightarrow \mathbb{R}^d$ is injective. Let $\Gamma$ be a group of congruence transformations preserving $\varPhi(X)$ and containing the period lattice. Clearly $\Gamma$ is a crystallographic group, so that we have an exact sequence: 
$$
0\longrightarrow \mathcal{L} \longrightarrow \Gamma\longrightarrow K\longrightarrow 1,
$$
where $\mathcal{L}(=\mathbb{R}^d\cap \Gamma)$ is a lattice containing the period lattice, and $K\subset O(d)$ is what is called the point group. Note that the isotropy group $\Gamma_{\bf x}$ (${\bf x}\in \mathbb{R}^d$) is identified with a subgroup of $K$ via the (injective) restriction of the homomorphism $\Gamma\longrightarrow K$ to $\Gamma_{\bf x}$.
Under the assumption on $\varPhi$, each $g\in \Gamma$ induces an automorphism of $X$, so $\Gamma$ is regarded as a subgroup of ${\rm Aut}(X)$.

\smallskip

\begin{thm}\label{thm:action}~ Suppose that 

\smallskip

{\rm (1)}~ the action of the point group $K$ on $\mathbb{R}^d$ is irreducible, and 
\smallskip

{\rm (2)}~ for any vertex $x\in V$, the fixed point set for $\Gamma_{\varPhi(x)}$-action on $\mathbb{R}^d$ is $\{{\bf 0}\}$, i.e., $\big\{{\bf x}\in \mathbb{R}^d|~g{\bf x}={\bf x}~(g\in \Gamma_{\varPhi(x)})\big\}=\{ {\bf 0}\}$. 

\smallskip

Then $\varPhi$ is a standard realization.   

\end{thm}
\noindent{\it Proof}.~ The group $K$ acts on $\big\{{\bf v}(e)\big\}_{e\in E_0}$ in a natural manner. Indeed, writing $\varPhi(\sigma x)=A(\sigma)\varPhi(x)+b(\sigma)$ ($\sigma\in \Gamma$), we have
$$
{\bf v}(\sigma e)=\varPhi\big(t(\sigma e)\big)-\varPhi\big(o(\sigma e)\big)=A(\sigma)\big[\varPhi\big(t(e)\big)-\varPhi\big(o(e)\big)\big]=A(\sigma){\bf v}(e).
$$
In view of Proposition \ref{thm:tightirre}, the assumption (1) assures us that $\big\{{\bf v}(e)\big\}_{e\in E_0}$ is tight. 
On the other hand, the vector 
$\displaystyle\sum_{e\in E_x}{\bf v}(e)$ is $\Gamma_{\varPhi(x)}$-invariant. Thus $\displaystyle\sum_{e\in E_x}{\bf v}(e)={\bf 0}$ by the assumption (2). This completes the proof. \hfill $\Box$

\medskip

As a corollary, we have

\begin{thm}
The $1$-skeleton of a Coxeter complex is a standard model.
\end{thm}

In this theorem, a Coxeter complex means a triangulation of $\mathbb{R}^d$ such that 
\smallskip

(a) if we denote by $\{\Delta_{\alpha}^d\}_{\alpha\in A}$ the set of $d$-dimensional simplices in this triangulation, then the action of  
the group $\Gamma$ generated by all reflections fixing facets of $\Delta_{\alpha}^d$ ($\alpha\in A$) preserves the triangulation, and 

\smallskip

(b) $\Gamma$ acts transitively on $\{\Delta_{\alpha}^d\}_{\alpha\in A}$.

\smallskip
For the proof, we use the fact that $\Gamma$ is a crystallographic group\footnote{Actually $\Gamma$ is isomorphic to a semidirect product of $\mathbb{Z}^d$ 
and the point group.}, which acts on the 1-skeleton as well, and that the point group for $\Gamma$ is the Weyl group associated with an irreducible root system. Clearly the fixed point set for $\Gamma_{{\bf x}}$-action is $\{{\bf 0}\}$ for every vertex ${\bf x}$, thereby the condition (2) being satisfied. Thus one may apply Theorem \ref{thm:action} to complete the proof. 
 
\medskip

A final remark is in order. The reader might wonder what is the practical use of standard crystal models. Straightforwardly speaking, 3D standard realizations  are purely mathematical outgrowths of logical reasoning. Even if a standard realization (or its deformation) looks realistic, it does not necessarily exist in nature; namely it is merely a model of a hypothetical crystal. Once we find a hypothetical crystal, however, a systematic prediction of its physical properties for appropriate atoms can be carried out by {\it first principles calculations} used in chemistry.
The prediction appealing to the computer power encourages (or discourages) material
scientists to synthesize the hypothetical crystals

\end{document}